\documentclass[10pt]{amsart}

\usepackage[margin=20pt,a4paper,top=3cm,bottom=4cm,left=2.5cm,right=2.5cm]{geometry}

\usepackage[english]{babel}
\usepackage{blindtext}
\usepackage{placeins}
\usepackage{eurosym}
\usepackage{subfig}
\usepackage{url}
\usepackage{fancyhdr}
\usepackage[utf8x]{inputenc}
\usepackage{indentfirst}
\usepackage{graphicx}
\usepackage{newlfont}
\usepackage{dsfont}
\usepackage{mathrsfs}
\usepackage{mathtools}
\mathtoolsset{showonlyrefs}
\usepackage{bbm}
\usepackage{amssymb}
\usepackage{amsmath}
\usepackage{amsfonts}
\usepackage{latexsym}
\usepackage{amsthm}
\usepackage{listings}
\usepackage{chngcntr}
\usepackage{imakeidx}

\def\de{\delta}
\def\R{\mathds{R}}

\def\p{\phi}
\def\wto{\rightharpoonup}

\def\po{\phi^\circ}
\def\pso{\psi^\circ}
\def\Div{\textup{div}}
\def\dist{\textup{dist}}

\def\e{\varepsilon}

\newcommand{\Om}{\Omega}

\renewcommand{\H}{\mathcal{H}}
\newcommand{\beq}{\begin{equation}}
\newcommand{\eeq}{\end{equation}}
\newcommand{\pa}{\partial}
\newcommand{\divr}{\mathfrak{Div}^r}

\newcommand{\vphi}{\varphi}
\newcommand{\osc}{\operatorname*{osc}}

\newcommand{\Lip}{\mathrm{Lip}}

\newcommand{\J}{J_r}
\renewcommand{\P}{\textnormal{Per}_r}
\newcommand{\ND}{\mathfrak{D}}
\newcommand{\td}{\tilde{d}}
\newcommand{\N}{\mathds N}

\newcommand{\res}
{\mathop{\hbox{\vrule height 7pt width .5pt depth 0pt
            \vrule height .5pt width 5pt depth 0pt}}\nolimits}

\theoremstyle{plain}

\newtheorem{theorem}{Theorem}[section]
\newtheorem{lemma}[theorem]{Lemma}
\newtheorem{proposition}[theorem]{Proposition}
\newtheorem{corollary}[theorem]{Corollary}
\theoremstyle{definition}
\begingroup
\newtheorem{remark}[theorem]{Remark}

\newtheorem{definition}[theorem]{Definition}
\endgroup
\numberwithin{equation}{section}

\usepackage{xcolor}
\usepackage{csquotes}

\usepackage{nomencl}
\makenomenclature

\makeindex
\usepackage{graphicx}

\title[Anisotropic Minkowski flows]{Anisotropic Minkowski flows: well-posedness and asymptotic behaviour}

\author{Giulia Anastasi}
\address{Department of Mathematical, Physical and Computer Sciences, University of Parma, Parco Area delle Scienze 53/A, 43124 Parma, Italy}
\email{giulia.anastasi@unipr.it}

\author{Filippo Cagnetti}
\address{Department of Mathematical, Physical and Computer Sciences, University of Parma, Parco Area delle Scienze 53/A, 43124 Parma, Italy}
\email{filippo.cagnetti@unipr.it}

\author{Massimiliano Morini}
\address{Department of Mathematical, Physical and Computer Sciences, University of Parma, Parco Area delle Scienze 53/A, 43124 Parma, Italy}
\email{massimiliano.morini@unipr.it}

\begin{document}

\begin{abstract}
We study anisotropic Minkowski-type curvature flows associated with a perimeter defined through oscillations on Wulff neighbourhoods. The anisotropy may be non-smooth, and in
particular crystalline, while the mobility may differ from the anisotropy, provided a natural compatibility condition is satisfied. We compute the first variation of the perimeter
and introduce a weak formulation in terms of distance functions and measure-valued nonlocal Cahn-Hoffman fields. By means of a nonlocal minimising movement scheme and suitable
estimates for the evolution of Wulff shapes, we prove existence of weak superflows and subflows starting from arbitrary closed sets. Under the additional assumption that the
initial level sets have compact boundary, a comparison principle yields existence and uniqueness of the level set flow. Finally, we prove that variational weak flows converge, up
to fattening, as $r\searrow 0$ to the weak flow of the local anisotropic, possibly crystalline, curvature motion. The main step in the proof consists in averaging the nonlocal
Cahn-Hoffman measures to obtain vector fields converging to a local Cahn-Hoffman field.
\end{abstract}

\maketitle

\section{Introduction}\label{sec:introduction}

Curvature motions are among the basic models describing the evolution of interfaces. Given a surface tension \(\p\) and a mobility \(\psi\), represented here by norms on \(\R^N\),
the corresponding geometric law is formally
\beq\label{intro local flow}
    V=-\psi(\nu)H^\p,
\eeq
where \(H^\p\) is the anisotropic mean curvature associated with the perimeter
\[
    P_\p(E):=\int_{\pa E}\p(\nu_E)\,d\H^{N-1}.
\]
For smooth and elliptic surface tensions, \eqref{intro local flow} has the usual differential geometric interpretation $H^\p=\Div_\tau\nabla\p(\nu)$, where $\Div_\tau$ denotes the
tangential divergence and $\nabla\p(\nu)$ is the so-called Cahn-Hoffman field. If \(\p\) is non-smooth, and in particular crystalline, facets may develop and the Cahn-Hoffman field
is no longer single valued. In this case a pointwise definition of curvature is no longer adequate. Following the approach in terms of signed distance functions introduced in
\cite{S93}, the distributional formulation developed in \cite{CMP17,CMNP19} replaces the pointwise curvature by the divergence of a suitable Cahn-Hoffman field taking values in
\(\pa\p(\nabla d)\). This formulation allows one to prove a comparison principle, convergence of the minimising movement scheme, and existence and uniqueness of the level set flow
also for crystalline surface tensions and anisotropic mobilities. An alternative approach, based on an extension of the viscosity theory to crystalline curvatures, has been
developed by Giga and Požár in recent years, see \cite{GP16,GP18,GP20} and the survey \cite{GP22}. We warn the reader that, with an abuse of language, we will refer to \eqref{intro local flow} as the local anisotropic flow although, strictly speaking, when $\phi$ is crystalline the curvature becomes a nonlocal object, 
as observed before.

In \cite{CMR26}, the aforementioned distributional approach has been developed in a nonlocal setting, namely for the fractional mean curvature flow and the isotropic Minkowski flow.

In this paper we extend the latter theory to anisotropic Minkowski flows.  
 For \(r>0\), let
\[
    W_r^\p(w):=\{x\in\R^N:\po(x-w)\le r\}
\]
be the Wulff-shape of radius \(r\) centred at \(w\), and set \(W_r^\p:=W_r^\p(0)\). We consider the anisotropic oscillation perimeter
\beq\label{intro perimeter}
    \P^\p(E) :=\frac1{2r}\int_{\R^N} \osc_{W_r^\p(w)}\chi_E\,dw.
\eeq
This functional admits an equivalent geometric interpretation as an anisotropic Minkowski pre-content. More precisely, for \(A\subset\R^N\) we denote by
\[
    (A)_r^\p := A+W_r^\p = \bigcup_{x\in A}\bigl(x+W_r^\p\bigr)
\]
the anisotropic tubular neighbourhood of \(A\) of radius \(r\). Then
\beq\label{intro minkowski precontent}
    \P^\p(E) = \frac1{2r} \mathcal L^N\!\left((\pa E^{(1)})_r^\p\right),
\eeq
where \(E^{(1)}\) denotes the set of points of density one of \(E\). Thus, \(\P^\p(E)\) is the normalized Lebesgue measure of an anisotropic tubular neighbourhood of the boundary
of the measure-theoretic representative of \(E\), with the geometry of the neighbourhood determined by the Wulff shape associated with \(\p\). The corresponding nonlocal total
variation is
\[
    \J^\p(u) :=\frac1{2r}\int_{\R^N} \osc_{W_r^\p(w)}u\,dw.
\]
The representation \eqref{intro minkowski precontent} shows how the interaction scale enters the definition of the energy. Indeed, irregularities of the boundary at scales much
smaller than \(r\) may give little or no additional contribution, while for smooth sets whose characteristic geometric scales are much larger than \(r\), the functional behaves
like a perimeter. This property is used in some variational models for image denoising, where one wants to remove small irregularities due to noise without affecting coherent fine
structures. It also appears in adversarial learning; see, for instance, \cite{BKLMP,CN17,BS24}.

The first variation of \eqref{intro perimeter} gives an anisotropic nonlocal curvature \(\mathcal K_r^\p\), and the geometric motion that we consider here is
\beq\label{curvature}
    V=-\psi(\nu)\mathcal K_r^\p.
\eeq
As explained below, such a curvature is affected by the cut locus and may fail to be defined on some smooth sets,
preventing a direct application of the general viscosity theory for translation-invariant nonlocal curvatures developed in \cite{CMP12,CMP15}.

In the isotropic case this difficulty was overcome in \cite{CMR26} by extending to the nonlocal setting the distributional formulation of \cite{CMP17,CMNP19}. More precisely, the
dual representation of the oscillation energy leads to the definition of suitable measure-valued nonlocal Cahn-Hoffman fields and of the corresponding nonlocal divergence. These
are then used to introduce a distributional notion of weak superflow and subflow in terms of distance functions, and to prove comparison, variational approximation, and existence
and uniqueness of the level set Minkowski flow. We also mention that closely related dual objects, together with nonlocal gradient and divergence operators and characterizations of
the subdifferential, have recently been introduced in the study of adversarial total variation in \cite{BS26}.

The aim of the present paper is twofold. First, for each fixed interaction scale, we extend the well-posedness theory of \cite{CMR26} to the anisotropic Minkowski motion
\eqref{curvature}. The anisotropy \(\p\) is arbitrary and may be crystalline, while the mobility \(\psi\) may differ from \(\p\) and is assumed to satisfy the natural
compatibility condition of \(\p\)-regularity introduced in \cite{CMNP19}. 

Second, we study the localisation of the Minkowski flow as \(r\downarrow0\), a problem not previously
addressed even for the isotropic Minkowski flow, and prove convergence to the local anisotropic, possibly crystalline, curvature motion \eqref{intro local flow} (with the same
surface tension \(\p\) and mobility \(\psi\)). At the level of the energies, this result is suggested by the \(\Gamma\)-convergence of anisotropic outer Minkowski contents to the
anisotropic perimeter established in \cite{CLL14}. Notice that an analogous result is proven in \cite{CDNP21} for the regularised version of the Minkowski flow introduced in
\cite{CMP12}, which can be studied by means of viscosity solutions. Related localisation results are also proved in \cite{CDGM26} for a fully discrete crystalline evolution, where
the spatial mesh determining the nonlocal interaction is sent to zero together with the time step, and in \cite{BLS24} for spatially inhomogeneous adversarial perimeters, where the
same small parameter determines both the interaction range and the time discretization. In the latter case the discrete evolutions converge to a weighted, locally isotropic
mean-curvature flow. In our case, we first construct the continuous weak evolution at every fixed scale \(r\), and then study the localisation limit of the corresponding flows.

\medskip

We now describe the results in more detail. In Section~\ref{sec:anisotropic-curvature} we identify the curvature appearing in \eqref{curvature}. Under regularity assumptions on
\(E\) and \(\p\), and outside a negligible exceptional set associated with the anisotropic cut locus, Proposition~\ref{curvature-prop} gives the first variation of \(\P^\p\). If
the inner and outer anisotropic normal segments of length \(r\) both realize the distance, then
\[
    \mathcal K_r^\p(E,x) =\frac{\det(I+rA_\p(x))-\det(I-rA_\p(x))}{2r}, \qquad A_\p(x):=\nabla^2\p(\nu_E(x))\nabla\nu_E(x).
\]
If one of the two normal segments reaches the cut locus before length \(r\), the corresponding Jacobian does not contribute to the curvature. At points with positive inner and
outer anisotropic reach, expansion of the determinants gives
\[
    \mathcal K_r^\p(E,x) =\operatorname{tr}A_\p(x)+O(r^2) =\Div_{\tau}\bigl(\nabla\p(\nu_E)\bigr)(x)+O(r^2) =H^\p(E,x)+O(r^2).
\]
Thus, also the first variation in the smooth case suggests the limiting law \eqref{intro local flow}.

In the first part of Section~\ref{sec:minkowski-flow} we study the generalised ROF problem
\beq\label{intro rof}
    \min_{v\in L^1_{\mathrm{loc}}(\R^N)} \left\{\J^\p(v)+\frac1{2h}\int_{\R^N}(v-g)^2\,dw\right\},
\eeq
or, more precisely, the associated Euler-Lagrange equation, which will be used to define the discrete geometric evolution. The Hilbert-space theory for data $g$ in \(\R+L^2\)
follows the strategy of \cite{CMR26}, whereas the argument for general Lipschitz data is different. By approximating an arbitrary Lipschitz datum by extensions which are constant
outside a compact set, without increasing its Lipschitz constant, Theorem~\ref{existence-uniqueness-Lip} yields a solution of the Euler-Lagrange system without any coercivity
assumption. The minimality property of the sublevel sets is formulated locally, with respect to compact perturbations, so that it remains meaningful even when the global nonlocal
perimeters are infinite. Proposition~\ref{comparison} establishes a one-sided comparison principle between coercive and possibly noncoercive data, as well as comparison with
suitable generalised Cahn-Hoffman supersolutions. As in \cite{CMR26}, uniqueness in the full Lipschitz class is obtained only when the datum is coercive.

In Section~\ref{min-mov-scheme-section} we introduce the nonlocal version of the classical minimising movement scheme \cite{ATW}, in the level set formulation of \cite{Chambolle}.
At each time step the signed distance from the current set is used as datum in the nonlocal ROF problem, and the set at the next step is defined as the zero sublevel set of a
solution. Thanks to the existence result for noncoercive data, the scheme can be defined for every closed initial set, thus extending the construction at a fixed scale in
\cite{CMR26}. Notice that for arbitrary closed sets we do not prove uniqueness of the solution at each step. Consequently, we do not establish a comparison principle in full
generality, although such a principle holds for compact initial data. The one-sided comparison of Proposition~\ref{comparison} is nevertheless sufficient to prove inclusion and
avoidance with respect to the uniquely determined evolutions of \(\psi\)-Wulff shapes. An important step in the proof is to estimate the evolution of these Wulff shapes. We first
control the discrete motion of \(\p\)-Wulff shapes; the \(\p\)-regularity of \(\psi\), in the sense of \cite{CMNP19}, then allows us to extend the estimates to \(\psi\)-Wulff
shapes. The resulting one-sided time estimates for the signed distances yield Helly-type compactness, and Theorem~\ref{convergence-min-mov-theo} shows that subsequential limits are
weak superflows and subflows of \eqref{curvature}. Section~\ref{sec:minkowski-flow} concludes with comparison and the well-posedness of the level set evolution for initial
functions whose sublevels have compact boundary, see Theorem~\ref{level-set-flow-theorem}.

Section~\ref{sec:asymptotic} is devoted to the asymptotic behaviour as $r\searrow 0$ and to the proof of the localisation theorem. The main difficulty is that the nonlocal
Cahn-Hoffman fields are measure-valued functions, while the local formulation involves vector fields. To overcome this difficulty, we associate with each nonlocal Cahn-Hoffman
field its first moment (see \eqref{z hat definitions}), which is a vector field. We then show that the nonlocal divergences converge to the ordinary divergence of a subsequential
weak-* limit $z$ of these first moments. One has also to prove that the limiting vector field $z$ is a local Cahn-Hoffman field for the limiting distance function; that is, it
satisfies
\[
    z\in\pa\p(\nabla d) \qquad\text{a.e. away from the evolving interface}.
\]
The latter condition is nonlinear. To prove it, we use a compensated compactness argument, which requires upper bounds, uniform in \(r\), for the nonlocal divergences on regions at
positive distance from the evolving interfaces. Such bounds follow from the minimising movement construction, and Theorem~\ref{asymptotic-behavior-theorem} consequently gives
subsequential convergence of variational weak superflows to weak superflows of \eqref{intro local flow}. However, in the Euclidean isotropic case,
Proposition~\ref{isotropic-divergence-bound} shows that the estimate holds for any family of weak superflows, even when they are not variational. Finally, under a nonfattening
assumption, the above results yield full convergence of variational nonlocal weak flows to the unique local weak flow.

\section{Preliminaries}\label{sec:preliminaries}

We collect here the notation and some basic facts from convex analysis and measure theory that will be used throughout the paper. We fix $N\in\N$. The spaces of Lipschitz and
compactly supported Lipschitz functions on $\R^N$ are denoted by $\Lip(\R^N)$ and $\Lip_c(\R^N)$, respectively, and
\[
\mathscr L^\pm:=\{g\in\Lip(\R^N): g(x)\to\pm\infty\ \text{as }|x|\to\infty\}.
\]
We write $\mathscr M(\R^N)$ for the Lebesgue measurable subsets of $\R^N$, $|A|=\mathcal L^N(A)$. For measurable functions we use the standard notation $f\vee g$, $f\wedge g$,
$f^+=f\vee0$, $f^-=(-f)^+$, and
\[
\R+L^2(\R^N):=\{g:\ g-M\in L^2(\R^N)\text{ for some }M\in\R\}.
\]
For a space-time set $E\subset\R^N\times[0,+\infty)$, $\mathrm{Int}(E)$ and $\overline E$ always refer to the relative topology of $\R^N\times[0,+\infty)$.

Let $\p$ be a norm on $\R^N$. Its polar norm is
\[
\po(\xi):=\sup_{\p(\nu)\le1}\nu\cdot\xi,
\]
and $(\po)^\circ=\p$; see, e.g., \cite{Rockafellar}. We shall repeatedly use the identity
\beq\label{sub gradient}
\pa\p(\nu)=\{\xi:\ \po(\xi)\le1,\ \xi\cdot\nu=\p(\nu)\}.
\eeq
In particular, if $\nu\neq0$ and $\xi\in\pa\p(\nu)$, then $\po(\xi)=1$.

For $r>0$ and $w\in\R^N$ we set
\[
W_r^\p(w):=\{\xi:\po(\xi-w)\le r\},\qquad W_r^\p:=W_r^\p(0),
\]
and call $W_r^\p(w)$ the $\p$-Wulff shape of radius $r$ centred at $w$. The distance induced by a norm $\p$ is
\beq\label{induced distance function}
\dist^\p(x,E):=\inf_{y\in E}\p(x-y),
\eeq
with the usual convention $\dist^\p(x,\emptyset)=+\infty$, and the signed distance is
\[
d_E^\p:=\dist^\p(\cdot,E)-\dist^\p(\cdot,E^c).
\]
Thus $\dist^\p(\cdot,E)=(d_E^\p)^+$ and $\dist^\p(\cdot,E^c)=(d_E^\p)^-$; moreover $\p(\nabla d_E^{\po})=\po(\nabla d_E^\p)=1$ a.e. away from $\pa E$.

We use Kuratowski convergence for closed sets. Namely, $E_n\overset{\mathscr K}{\longrightarrow}E$ if every limit of points $x_n\in E_n$ belongs to $E$ and every $x\in E$ is the
limit of some $x_n\in E_n$. Equivalently, for any norm $\p$,
\[
\dist^\p(\cdot,E_n)\longrightarrow\dist^\p(\cdot,E) \quad\text{locally uniformly in }\R^N.
\]
In particular every sequence of closed sets has a Kuratowski convergent subsequence.

Let $\mathcal M(A)$ and $\mathcal M^+(A)$ denote the signed and positive Radon measures on a Borel set $A$, respectively. For $r>0$ set
\beq\label{Wulff product shape}
C_r^\p:=W_r^\p\times W_r^\p.
\eeq
The anisotropic $r$-neighbourhood of $E$ is $(E)_r^\p:=\bigcup_{w\in E}W_r^\p(w)$. Since $C_r^\p$ is compact, $(C^0(C_r^\p))'=\mathcal M(C_r^\p)$. We shall use the standard
identification (see \cite[Definition~2.111]{FonsecaLeoniBook})
\[
\big(L^1(\R^N;C^0(C_r^\p))\big)'=L^\infty_{\mathbf w}(\R^N;\mathcal M(C_r^\p)),
\]
where $z\in L^\infty_{\mathbf w}$ means that $w\mapsto\int f\,dz_w$ is measurable for every $f\in C^0(C_r^\p)$ and $\operatorname*{ess\,sup}_w|z_w|(C_r^\p)<\infty$, and $z_w$ stands for the evaluation of $z$ at $w$.

The anisotropic Minkowski-type perimeter and the corresponding nonlocal total variation are
\beq\label{perimeter def}
\mathrm{Per}_r^\p(E):=\frac1{2r}\int_{\R^N}\osc_{W_r^\p(w)}\chi_E\,dw,
\eeq
\[
\J^\p(v):=\frac1{2r}\int_{\R^N}\osc_{W_r^\p(w)}v\,dw.
\]
They are the Wulff-neighbourhood versions of the functionals considered in \cite{CMR26}; see also \cite{CLL14}. We will use the generalised coarea formula
\beq\label{coarea formula}
\J^\p(v)=\int_{-\infty}^{+\infty}\mathrm{Per}_r^\p(\{v>t\})\,dt =\int_{-\infty}^{+\infty}\mathrm{Per}_r^\p(\{v\le t\})\,dt,
\eeq
as well as translation invariance, one-homogeneity, subadditivity, truncation monotonicity, and the submodularity identities
\beq\label{sub modularity global}
\mathrm{Per}_r^\p(E\cap F)+\mathrm{Per}_r^\p(E\cup F) \le \mathrm{Per}_r^\p(E)+\mathrm{Per}_r^\p(F),
\eeq
\[
\J^\p(u\vee v)+\J^\p(u\wedge v)\le \J^\p(u)+\J^\p(v).
\]
These properties follow exactly as in \cite[Section~2]{CMR26}. The functional $\J^\p$ is convex and lower semicontinuous under strong, hence also weak, $L^p$ convergence for
$1\le p<\infty$.

We conclude by recalling the distributional formulation, in terms of distance functions, of weak flows for the local anisotropic, possibly crystalline, motion
\beq\label{curvature00}
V=-\psi(\nu)H^\p.
\eeq
The following definitions are taken from \cite{CMP17} and, for general mobilities, from \cite{CMNP19}.

\begin{definition}\label{classic-superflow-def}
Let $E^0\subset\R^N$ be closed and let $E\subset\R^N\times[0,+\infty)$ be closed, with $E(t):=\{x:(x,t)\in E\}$. We say that $E$ is a weak superflow of \eqref{curvature00} with
initial datum $E^0$ if the following conditions hold:
\begin{itemize}
\item[(a)] $E(0)\subset E^0$;
\item[(b)] $E(s)\overset{\mathscr K}{\longrightarrow}E(t)$ as $s\nearrow t$ for every $t>0$;
\item[(c)] after extinction the flow remains empty;
\item[(d)] setting $T^*:=\inf\{t>0:E(s)=\emptyset\text{ for }s\ge t\}$ and $d(x,t):=\dist^{\pso}(x,E(t))$, there exists $z\in L^\infty(\R^N\times(0,T^*);\R^N)$ such that
\[
\pa_t d\ge\Div z\quad\text{in }(\R^N\times(0,T^*))\setminus E,
\]
in the distributional sense, $z\in\pa\p(\nabla d)$ a.e., and $(\Div z)^+$ is bounded on every set $\{d\ge\de\}$, $\de>0$.
\end{itemize}
An open set $A$ is a weak subflow if its complement is a weak superflow with the complementary initial datum. If $E^0=\overline{\mathrm{Int}(E^0)}$, a weak flow is a closed set
$E=\overline{\mathrm{Int}(E)}$ which is a weak superflow and whose interior is a weak subflow.
\end{definition}

\begin{definition}[Level set solutions]\label{sup-sub-levle-set-def-classic}
Let $u^0$ be uniformly continuous. A lower semicontinuous $u$ is a level set supersolution of \eqref{curvature00} if $u(\cdot,0)\ge u^0$ and, for a.e. $\lambda$, $\{u\le\lambda\}$
is a weak superflow with initial datum $\{u^0\le\lambda\}$. Upper semicontinuous subsolutions are defined by changing sign, and a continuous function which is both a sub- and a
supersolution is called a level set solution.
\end{definition}

\section{The anisotropic nonlocal curvature}\label{sec:anisotropic-curvature}
In this section we compute the first variation of \eqref{perimeter def}. Throughout the section $\p$ is a $C^2$ norm on $\R^N \setminus \{ 0 \}$ and $E$ is a $C^2$ set with compact boundary. We
denote by $\nu_E$ the Euclidean outer unit normal and extend it smoothly to a tubular neighbourhood of $\pa E$ by $\nu_E=\nabla d_E$.

For $x\in\pa E$ set
\[
\mathcal K_r^\p(E,x):=\mathcal K_r^{\p,+}(E,x)+\mathcal K_r^{\p,-}(E,x),
\]
where
\[
\mathcal K_r^{\p,+}(E,x):=
\begin{cases}
\displaystyle \frac1{2r}\det\!\bigl(I+r\nabla^2\p(\nu_E(x))\nabla\nu_E(x)\bigr), &\text{if }\dist^{\po}\!\bigl(x+r\nabla\p(\nu_E(x)),E\bigr)=r,\\[1ex]
0,&\text{otherwise},
\end{cases}
\]
and
\[
\mathcal K_r^{\p,-}(E,x):=
\begin{cases}
\displaystyle -\frac1{2r}\det\!\bigl(I-r\nabla^2\p(\nu_E(x))\nabla\nu_E(x)\bigr), &\text{if }\dist^{\po}\!\bigl(x-r\nabla\p(\nu_E(x)),E^c\bigr)=r,\\[1ex]
0,&\text{otherwise}.
\end{cases}
\]
Thus the two terms describe, respectively, the first variation of the outer and inner anisotropic strips of width $r$.

For the outer strip we decompose $\pa E=A_r^+\cup B_r^+\cup\mathcal N_r^+$, where
\[
A_r^+:=\Bigl\{x\in\pa E:\ \dist^{\po}\!\bigl(x+t\nabla\p(\nu_E(x)),E\bigr)=t \text{ for some }t>r\Bigr\},
\]
\[
B_r^+:=\Bigl\{x\in\pa E:\ \dist^{\po}\!\bigl(x+r\nabla\p(\nu_E(x)),E\bigr)<r\Bigr\}, \qquad \mathcal N_r^+:=\pa E\setminus(A_r^+\cup B_r^+).
\]
The sets $A_r^-$, $B_r^-$ and $\mathcal N_r^-$ are defined analogously by looking in the direction $-\nabla\p(\nu_E)$ and replacing $E$ by $E^c$ in the distance condition. We set
$\mathcal N_r:=\mathcal N_r^+\cup\mathcal N_r^-$.

\begin{proposition}\label{curvature-prop}
Let $E$ be as above and assume $\H^{N-1}(\mathcal N_r)=0$. For every $\vphi\in C^2(\pa E)$,
\[
\left.\frac d{d\e}\mathrm{Per}_r^\p(\Phi_\e(E))\right|_{\e=0} = \int_{\pa E}\mathcal K_r^\p(E,x)\,\vphi(x)\p(\nu_E(x))\,d\H^{N-1}(x),
\]
where $\Phi_\e$ is any smooth family of diffeomorphisms satisfying
\begin{equation}\label{displacement-def}
\Phi_\e(x)=x+\e\vphi(x)\nabla\p(\nu_E(x)) \qquad x\in\pa E.
\end{equation}
\end{proposition}

We first prove the following expansion of the anisotropic normal map. We extend $\vphi$ smoothly to a neighbourhood of $\pa E$, keeping it constant along the Euclidean normal
lines, so that $\nabla\vphi=\nabla_\tau\vphi$ on $\pa E$. Throughout the lemma and the proof of Proposition~\ref{curvature-prop}, all remainders denoted by $O(\e^2)$ or $o(\e)$ are
uniform with respect to the variables ranging in the compact sets under consideration. Thus, for instance, $R(x,t,\e)=o(\e)$ uniformly for $x\in K$ and $t\in I$ means
\[
\sup_{x\in K,\ t\in I}\frac{|R(x,t,\e)|}{|\e|}\longrightarrow0 \qquad\text{as }\e\to0.
\]

\begin{lemma}\label{curvature-lemma}
Let $\vphi\in C^2(\pa E)$ and define
\[
G_\e(x,t):=\Phi_\e(x)+t\nabla\p\!\left(\nu_{\Phi_\e(E)}(\Phi_\e(x))\right).
\]
Then, uniformly for $x\in\pa E$ and $|t|\le r$,
\beq\label{0-0}
G_\e(x,t) =x+(t+\e\vphi(x))\nabla\p(\nu_E(x)) -\e t\,\p(\nu_E(x))\nabla^2\p(\nu_E(x))\nabla_\tau\vphi(x) +o(\e).
\eeq
\end{lemma}

\begin{proof}
Set
\[
A(x):=\nabla\p(\nu_E(x))\otimes\nabla_\tau\vphi(x) +\vphi(x)\nabla^2\p(\nu_E(x))\nabla\nu_E(x).
\]
On $\pa E$ we have $D\Phi_\e=I+\e A$. Hence
\[
(D\Phi_\e)^{-T}=I-\e A^T+O(\e^2)
\]
uniformly on $\pa E$. By the standard transformation formula for the normal (see, e.g., \cite[Eq.~(2.6)]{CMM08}),
\[
\nu_{\Phi_\e(E)}(\Phi_\e(x)) =\frac{(D\Phi_\e)^{-T}\nu_E(x)}{|(D\Phi_\e)^{-T}\nu_E(x)|}.
\]
Euler's identity and the one-homogeneity of $\p$ give
\[
\nabla\p(\nu_E)\cdot\nu_E=\p(\nu_E), \qquad \nabla^2\p(\nu_E)\nu_E=0.
\]
Since $\nabla_\tau\vphi\perp\nu_E$, we obtain
\[
A^T\nu_E=\p(\nu_E)\nabla_\tau\vphi.
\]
Consequently,
\[
(D\Phi_\e)^{-T}\nu_E =\nu_E-\e\p(\nu_E)\nabla_\tau\vphi+O(\e^2).
\]
Notice that the normalization in the formula for the transformed normal does not change the first-order term. Indeed, since $\nu_E\cdot\nabla_\tau\vphi=0$,
\[
\bigl|(D\Phi_\e)^{-T}\nu_E\bigr|^2 =\bigl|\nu_E-\e\p(\nu_E)\nabla_\tau\vphi+O(\e^2)\bigr|^2 =1+O(\e^2),
\]
so that
\[
\bigl|(D\Phi_\e)^{-T}\nu_E\bigr|^{-1}=1+O(\e^2).
\]
Therefore
\beq\label{normal variation}
\nu_{\Phi_\e(E)}(\Phi_\e(x)) =\nu_E(x)-\e\p(\nu_E(x))\nabla_\tau\vphi(x)+O(\e^2).
\eeq
Because $\p\in C^2$, Taylor expansion of the $C^1$ map $\nabla\p$ yields, uniformly on the compact set $\nu_E(\pa E)$,
\beq\label{anisotropic normal variation}
\nabla\p\!\left(\nu_{\Phi_\e(E)}(\Phi_\e(x))\right) =\nabla\p(\nu_E(x)) -\e\p(\nu_E(x))\nabla^2\p(\nu_E(x))\nabla_\tau\vphi(x) +o(\e).
\eeq
Substituting \eqref{anisotropic normal variation} in the definition of $G_\e$ proves \eqref{0-0}.
\end{proof}

\begin{proof}[Proof of Proposition~\ref{curvature-prop}]
For a set $F$ define the outer and inner contributions relative to $E$ by
\[
\mathrm{Per}_{r}^{\p,+}(F\mid E):=\frac1{2r}\int_{E^c}\osc_{W_r^\p(w)}\chi_F\,dw, \qquad \mathrm{Per}_{r}^{\p,-}(F\mid E):=\frac1{2r}\int_E\osc_{W_r^\p(w)}\chi_F\,dw.
\]
At $F=E$ their sum is $\mathrm{Per}_r^\p(E)$. We prove the formula for the outer contribution; the inner one will then follow by passing to the complement.

\smallskip
\noindent\textbf{Step 1: perturbations supported in $A_r^+$.}
Assume $\operatorname{supp}\vphi\Subset A_r^+$. By compactness there exists $\eta>0$ such that the points given by the anisotropic normal map
\[
F(x,s):=x+s\nabla\p(\nu_E(x))
\]
have a unique projection onto $E$ with respect to the anisotropic distance, for $x\in\operatorname{supp}\vphi$ and $0\le s\le r+\eta$. In particular,
\beq\label{anisotropic distance tube}
d_E^{\po}(F(x,s))=s, \qquad \nabla d_E^{\po}(F(x,s))=\frac{\nu_E(x)}{\p(\nu_E(x))}.
\eeq
We now use the expansion of Lemma~\ref{curvature-lemma} to compute the displacement of the outer boundary of the anisotropic neighbourhood along the anisotropic normal directions.
Set
\[
\de_\e(x,t):=G_\e(x,t)-F(x,t).
\]
By \eqref{0-0},
\[
\de_\e(x,t) =\e\vphi(x)\nabla\p(\nu_E(x)) -\e t\,\p(\nu_E(x))\nabla^2\p(\nu_E(x))\nabla_\tau\vphi(x) +o(\e).
\]
On the compact tubular region under consideration the function $d_E^{\po}$ is $C^1$ with uniformly continuous gradient. Taylor expansion at $F(x,t)$, together with
\eqref{anisotropic distance tube}, therefore gives
\[
\begin{split}
d_E^{\po}(G_\e(x,t)) &=d_E^{\po}(F(x,t)) +\nabla d_E^{\po}(F(x,t))\cdot\de_\e(x,t) +o(\e)\\
&=t+\frac{\nu_E(x)}{\p(\nu_E(x))}\cdot\de_\e(x,t) +o(\e).
\end{split}
\]
We compute the two first-order terms in the last scalar product separately. By Euler's identity,
\[
\frac{\nu_E}{\p(\nu_E)}\cdot \bigl(\vphi\nabla\p(\nu_E)\bigr)=\vphi,
\]
whereas
\[
\frac{\nu_E}{\p(\nu_E)}\cdot \bigl(-t\p(\nu_E)\nabla^2\p(\nu_E)\nabla_\tau\vphi\bigr) =-t\,\nabla^2\p(\nu_E)\nu_E\cdot\nabla_\tau\vphi=0.
\]
Hence
\beq\label{0-9}
d_E^{\po}(G_\e(x,t)) =t+\e\vphi(x)+o(\e).
\eeq
We now describe the outer boundary of the anisotropic $r$-neighbourhood of the perturbed set along the \emph{fixed} anisotropic normal directions of $E$. Set $E_\e:=\Phi_\e(E)$
and, for $(x,s)$ in a fixed neighbourhood of $\operatorname{supp}\vphi\times\{r\}$, define
\beq\label{definition Hepsilon}
H_\e(x,s):=d_{E_\e}^{\po}\bigl(F(x,s)\bigr)-r.
\eeq
We first justify the regular dependence of these distance functions on the parameter $\e$. To this aim consider the map
\[
\mathcal G(\e,y,t):=\bigl(\e,G_\e(y,t)\bigr).
\]
It is of class $C^1$ in a neighbourhood of the compact set under consideration. Indeed, $(\e,y)\mapsto\Phi_\e(y)$ is $C^1$, and
\[
\nu_{E_\e}(\Phi_\e(y)) =\frac{(D\Phi_\e(y))^{-T}\nu_E(y)}{|(D\Phi_\e(y))^{-T}\nu_E(y)|};
\]
since $\p\in C^2$, it follows that $(\e,y)\mapsto\nabla\p(\nu_{E_\e}(\Phi_\e(y)))$ is $C^1$ as well. At $\e=0$ we have $G_0=F$, and the differential of $\mathcal G$ has the block
form
\[
D\mathcal G(0,y,t)=
\begin{pmatrix}
1 & 0\\
\left.\pa_\e G_\e(y,t)\right|_{\e=0} & D_{(y,t)}F(y,t)
\end{pmatrix}.
\]
The lower-right block is invertible throughout the tubular region under consideration. Hence the inverse function theorem, applied to $\mathcal G$, yields, after possibly shrinking
the neighbourhood and using compactness, $C^1$ inverse coordinates of the form
\[
(\e,w)\longmapsto \bigl(\e,Y(\e,w),T(\e,w)\bigr), \qquad w=G_\e\bigl(Y(\e,w),T(\e,w)\bigr).
\]
In the part of the tubular neighbourhood outside the set, the second coordinate is the anisotropic distance from $E_\e$, namely
\[
T(\e,w)=d_{E_\e}^{\po}(w).
\]
Consequently $(\e,w)\mapsto d_{E_\e}^{\po}(w)$, and therefore $(\e,x,s)\mapsto H_\e(x,s)$, are $C^1$ in the region under consideration. At $\e=0$ we have
\[
H_0(x,r)=0
\]
and, by \eqref{anisotropic distance tube} and Euler's identity,
\beq\label{Hs at zero}
\pa_s H_0(x,r) =\nabla d_E^{\po}(F(x,r))\cdot\pa_sF(x,r) =\frac{\nu_E}{\p(\nu_E)}\cdot\nabla\p(\nu_E)=1.
\eeq
Therefore, by the implicit function theorem, for $|\e|$ small there exists a unique $C^1$ function $s_\e(x)$, defined for $x$ in a fixed neighbourhood of
$\operatorname{supp}\vphi$, such that
\beq\label{definition sepsilon}
H_\e\bigl(x,s_\e(x)\bigr)=0, \qquad s_0(x)=r.
\eeq
Equivalently, $F(x,s_\e(x))$ is the point where the outer boundary of the anisotropic $r$-neighbourhood of $E_\e$ meets the fixed ray through $x$, namely
\[
\{x+s\nabla\p(\nu_E(x)):s\ge0\}.
\]

It remains to compute the first-order displacement of this intersection. Fix $x$ and use the inverse coordinates just constructed at the fixed point $F(x,r)$: set
\[
Y_\e(x):=Y\bigl(\e,F(x,r)\bigr), \qquad T_\e(x):=T\bigl(\e,F(x,r)\bigr) =d_{E_\e}^{\po}(F(x,r)).
\]
Thus,
\beq\label{F equal G_e}
F(x,r)=G_\e\bigl(Y_\e(x),T_\e(x)\bigr), \qquad Y_0(x)=x, \qquad T_0(x)=r.
\eeq
We denote differentiation with respect to $\e$ at $\e=0$ by a dot; in particular,
\[
\dot T_0(x):=\left.\pa_\e T_\e(x)\right|_{\e=0} =\left.\pa_\e H_\e(x,r)\right|_{\e=0}.
\]
Differentiating \eqref{F equal G_e} at $\e=0$ and using the preceding identity, we obtain
\[
0=\left.\pa_\e G_\e(x,r)\right|_{\e=0} +D_xF(x,r)[\dot Y_0(x)] +\nabla\p(\nu_E(x))\,\dot T_0(x).
\]
Taking the scalar product with
\[
\nabla d_E^{\po}(F(x,r))=\frac{\nu_E(x)}{\p(\nu_E(x))}
\]
eliminates the middle term, since $d_E^{\po}(F(y,r))=r$ and hence
\[
\nabla d_E^{\po}(F(x,r))\cdot D_xF(x,r)[\tau]=0 \qquad\text{for every }\tau\in T_x\pa E.
\]
Moreover, by Lemma~\ref{curvature-lemma},
\[
\left.\pa_\e G_\e(x,r)\right|_{\e=0} =\vphi(x)\nabla\p(\nu_E(x)) -r\p(\nu_E(x))\nabla^2\p(\nu_E(x))\nabla_\tau\vphi(x),
\]
and hence, exactly as in the proof of \eqref{0-9},
\[
\frac{\nu_E(x)}{\p(\nu_E(x))}\cdot \left.\pa_\e G_\e(x,r)\right|_{\e=0} =\vphi(x).
\]
Using also \eqref{Hs at zero}, we conclude that
\[
\dot T_0(x)=-\vphi(x), \qquad\text{that is,}\qquad \left.\pa_\e H_\e(x,r)\right|_{\e=0}=-\vphi(x).
\]
Differentiating \eqref{definition sepsilon} at $\e=0$ now yields
\[
0=-\vphi(x)+\pa_sH_0(x,r)\,\left.\pa_\e s_\e(x)\right|_{\e=0} =-\vphi(x)+\left.\pa_\e s_\e(x)\right|_{\e=0}.
\]
Therefore
\beq\label{outer offset graph}
s_\e(x)=r+\e\vphi(x)+o(\e),
\eeq
uniformly on $\operatorname{supp}\vphi$, the uniformity following from the $C^1$ dependence above and compactness. This is the only description of the perturbed outer boundary
needed in the volume computation below.

We next compute the volume element in the fixed tubular coordinates. Let $\tau_1,\ldots,\tau_{N-1}$ be an oriented orthonormal basis of $T_x\pa E$. For every tangential vector
$\tau$,
\[
D_xF(x,s)[\tau] =\tau+s\nabla^2\p(\nu_E(x))\nabla\nu_E(x)[\tau], \qquad \pa_sF(x,s)=\nabla\p(\nu_E(x)).
\]
Since $\nabla^2\p(\nu_E)\nu_E=0$ and $\nabla^2\p(\nu_E)$ is symmetric,
\[
\nu_E(x)\cdot D_xF(x,s)[\tau]=0;
\]
hence the first $N-1$ vectors remain in $T_x\pa E$. Decomposing $\nabla\p(\nu_E)$ into tangential and normal components, its tangential component gives no contribution to the
determinant, and Euler's identity yields
\[
\nabla\p(\nu_E)\cdot\nu_E=\p(\nu_E).
\]
Consequently the full Jacobian of $(x,s)\mapsto F(x,s)$ is
\beq\label{full tube jacobian}
J_F(x,s) =\p(\nu_E(x))\, \det_T\!\bigl(I+s\nabla^2\p(\nu_E(x))\nabla\nu_E(x)\bigr).
\eeq
Here $\det_T$ denotes the determinant on $T_x\pa E$. Since the extension of $\nu_E$ is constant along Euclidean normal lines, $\nabla\nu_E\,\nu_E=0$. Moreover the operator inside
the determinant preserves $T_x\pa E$ by the orthogonality relation proved above, and acts as the identity in the normal direction. Therefore its full determinant in $\R^N$
coincides with its tangential determinant, and we may write
\beq\label{tube jacobian}
J_F(x,s) =\p(\nu_E(x))\, \det\!\bigl(I+s\nabla^2\p(\nu_E(x))\nabla\nu_E(x)\bigr).
\eeq
The determinant is positive on the region under consideration: it equals $1$ at $s=0$ and cannot vanish before the anisotropic normal parametrization loses local invertibility,
while the latter remains regular up to $r+\eta$.

We are now in a position to compute the variation of the outer contribution. For $F=E$ and $w\in E^c$ in the above tubular region, $\osc_{W_r^\p(w)}\chi_E=1$ precisely for $0<s<r$
(up to the boundary levels, which have zero measure). For $E_\e:=\Phi_\e(E)$ and $|\e|$ small, the inner boundary $s=0$ of the integration region $E^c$ remains fixed, while the
outer boundary moves to $s=s_\e(x)$. Outside a fixed neighbourhood of $\operatorname{supp}\vphi$ the integrand is unchanged. Therefore the difference between the two outer
contributions is, up to $o(\e)$, the signed volume of the thin strip between $s=r$ and $s=s_\e(x)$. Using \eqref{tube jacobian} and the area formula (equivalently, coarea together
with \eqref{anisotropic distance tube}; see \cite[Theorems~2.91 and 2.93]{AFP}), we obtain
\beq\label{incremental ratio}
\begin{split}
&\frac{\mathrm{Per}_{r}^{\p,+}(\Phi_\e(E)\mid E)-\mathrm{Per}_{r}^{\p,+}(E\mid E)}{\e}\\
&\quad=\frac1{2r}\int_{\pa E} \frac1\e\int_r^{r+\e\vphi(x)+o(\e)} \det\!\bigl(I+s\nabla^2\p(\nu_E(x))\nabla\nu_E(x)\bigr) \p(\nu_E(x))\,ds\,d\H^{N-1}(x)+o(1).
\end{split}
\eeq
Here the integrand is unchanged outside a fixed neighbourhood of $\operatorname{supp}\vphi$. Letting $\e\to0$ gives
\beq\label{0-2}
\left.\frac d{d\e}\mathrm{Per}_{r}^{\p,+}(\Phi_\e(E)\mid E)\right|_{\e=0} =\int_{\pa E}\mathcal K_r^{\p,+}(E,x)\vphi(x)\p(\nu_E(x))\,d\H^{N-1}(x)
\eeq
whenever $\operatorname{supp}\vphi\Subset A_r^+$. The argument applies to both signs of $\e$; equivalently, the negative one-sided derivative follows by replacing $\vphi$ with
$-\vphi$.

\smallskip
\noindent\textbf{Step 2: perturbations supported in $B_r^+$.}
Assume $\operatorname{supp}\vphi\Subset B_r^+$. By compactness there exists $\eta>0$ such that
\[
\dist^{\po}\!\bigl(x+r\nabla\p(\nu_E(x)),E\bigr)\le r-\eta \qquad\text{for }x\in\operatorname{supp}\vphi.
\]
The strict inequality is preserved under sufficiently small $C^1$ perturbations of the boundary. Consequently, the corresponding points along the normal directions remain in the
interior of the $r$-neighbourhood, and this part of the boundary gives no first-order contribution to the variation of the outer strip. Hence
\[
\left.\frac d{d\e}\mathrm{Per}_{r}^{\p,+}(\Phi_\e(E)\mid E)\right|_{\e=0}=0
\]
for such $\vphi$, in agreement with $\mathcal K_r^{\p,+}=0$ on $B_r^+$.

\smallskip
\noindent\textbf{Step 3: localisation near $\mathcal N_r^+$.}
Let $\de>0$. Since $\H^{N-1}(\mathcal N_r^+)=0$, choose an open neighborhood $U_r^+$ of $\mathcal N_r^+$ in $\pa E$ such that
\[
\H^{N-1}(U_r^+)\le\de.
\]
Take a smooth partition of unity $\eta_1+\eta_2+\eta_3=1$ on $\pa E$ subordinated to $A_r^+$, $B_r^+$ and $U_r^+$, and set $\vphi_i:=\eta_i\vphi$. Moreover, we denote by $\Phi_\e^{\vphi_1+\vphi_2}$ the diffeomorphism defined in \eqref{displacement-def} with $\vphi_1+\vphi_2$ in place of $\vphi$. The area formula in a fixed
tubular neighbourhood, together with the uniform $C^1$ bound on the deformations, gives
\beq\label{localization estimate}
\left| \bigl[\mathrm{Per}_{r}^{\p,+}(\Phi_\e(E)\mid E)-\mathrm{Per}_{r}^{\p,+}(E\mid E)\bigr] - \bigl[\mathrm{Per}_{r}^{\p,+}(\Phi_\e^{\vphi_1+\vphi_2}(E)\mid
E)-\mathrm{Per}_{r}^{\p,+}(E\mid E)\bigr] \right| \le C|\e|\de,
\eeq
where $C$ depends on $E$, $r$, $\p$ and $\|\vphi\|_{C^1}$, but not on $\de$ or $\e$. Steps 1 and 2 applied to $\vphi_1$ and $\vphi_2$, together with \eqref{localization estimate},
imply
\[
\limsup_{\e\to0}\left| \frac{\mathrm{Per}_{r}^{\p,+}(\Phi_\e(E)\mid E)-\mathrm{Per}_{r}^{\p,+}(E\mid E)}{\e} - \int_{\pa E}\mathcal
K_r^{\p,+}(E,x)\vphi(x)\p(\nu_E(x))\,d\H^{N-1}(x) \right|\le C\de.
\]
Letting $\de\downarrow0$ proves the outer variation formula.

\smallskip
\noindent\textbf{Step 4: the inner contribution.}
Observe that
\[
\mathrm{Per}_{r}^{\p,-}(F\mid E) =\mathrm{Per}_{r}^{\p,+}(F^c\mid E^c).
\]
Since $\Phi_\e(E)^c=\Phi_\e(E^c)$, we apply the outer formula to $E^c$. On $\pa E$,
\[
\nu_{E^c}=-\nu_E, \qquad \nabla\p(-\nu_E)=-\nabla\p(\nu_E), \qquad \nabla^2\p(-\nu_E)=\nabla^2\p(\nu_E), \qquad \nabla\nu_{E^c}=-\nabla\nu_E.
\]
Moreover the same geometric deformation is written relative to the outer anisotropic normal of $E^c$ with scalar speed $-\vphi$. Therefore
\[
\left.\frac d{d\e}\mathrm{Per}_{r}^{\p,-}(\Phi_\e(E)\mid E)\right|_{\e=0} =\int_{\pa E}\mathcal K_r^{\p,-}(E,x)\vphi(x)\p(\nu_E(x))\,d\H^{N-1}(x).
\]
Adding this identity to the formula for the outer variation concludes the proof.
\end{proof}

\section{The Minkowski-type mean curvature flow}\label{sec:minkowski-flow}
\subsection{The generalised ROF problem.} The aim of this section is to study, for $h >0$ fixed, the minimum problem
\beq\label{minimum problem}
    \min_{ v \in L^1_{\mathrm{loc}}(\R^N)} \Big( \J^\p(v) +\frac{1}{2h} \int_{\R^N} (v-g)^2 \, dw \Big),
\eeq
and the associated Euler-Lagrange equation. We start by considering the case $g \in \R + L^2(\R^N)$.
\begin{proposition}\label{L2-existence-uniqueness}
    Let $h>0$. Then, for every $g \in \R + L^2(\R^N)$ there exists a unique solution $u_h^g$ to the problem \eqref{minimum problem}, with $u_h^g \in \R + L^2(\R^N)$. Moreover, if
    $g_1, g_2 \in \R + L^2(\R^N)$ and $g_1 \leq g_2$ a.e. in $\R^N$, then $u_h^{g_1} \leq u_h^{g_2}$ a.e. in $\R^N$. Finally, if $g \in \Lip(\R^N) \cap ( \R + L^2(\R^N))$, then
    $u_h^g \in \Lip(\R^N) \cap ( \R + L^2(\R^N))$ with $\Lip(u_h^g) \leq \Lip(g)$, and there exists $z \in L^{\infty}_{\mathbf{w}}(\R^N; \mathcal{M}^+(C_r^\p))$ such that
    \begin{itemize}
        \item[(a)] for a.e. $w \in \R^N$ one has $z_w(C_r^\p) \leq 1$ and
        \beq\label{E L equation L2}
            \int_{C_r^\p} (u_h^g(w+x)-u_h^g(w+y)) \, dz_w(x,y) = \mathop{\mathrm{osc}}_{W_r^{\p}(w)} u_h^g;
        \eeq
        \item[(b)] for every $\vphi \in C^0(\R^N) \cap L^2(\R^N)$, with $\J^\p(\vphi) < + \infty$
        \beq\label{Cahn Hoffman L2}
        \frac{h}{2r} \int_{\R^N} \int_{C_r^\p} (\vphi(w+x)-\vphi(w+y)) \, dz_w(x,y) = -\int_{\R^N} \vphi(u_h^g-g) \, dw,
        \eeq
    \end{itemize}
    where $z_w$ stands for the evaluation of $z$ at the point $w$, and $dz_w(x,y)$ for the integration in $(x,y)$ with respect to the measure $z_w$.
\end{proposition}
\begin{proof}
    We can argue exactly as in \cite[Proposition 3.1]{CMR26}.
\end{proof}

Following \cite[Section~3.1]{CMR26}, we introduce the anisotropic nonlocal divergence of measure-valued Cahn-Hoffman fields, replacing Euclidean balls by $\p$-Wulff shapes. For
$z\in L^{\infty}_{\mathbf w}(\R^N;\mathcal M(C_r^\p))$, the measure $\divr_\p z\in\mathcal M(\R^N)$ is defined by
\beq
    \int_{\R^N} \vphi(w) \, d(\divr_\p z)(w) \coloneqq - \frac{1}{2r} \int_{\R^N} \int_{C_r^\p} (\vphi(w+x)-\vphi(w+y)) \, dz_w(x,y)dw
\eeq
for all $\vphi\in C^0_c(\R^N)$. We recall that $\mathcal{M}(\R^N)$ denotes the space of Radon measures on $\R^N$.

\begin{definition}\label{generalized-Cahn-Hoffman-definition}
Let $u\in\Lip(\R^N)$. A field $z\in L^\infty_{\mathbf w}(\R^N;\mathcal M^+(C_r^\p))$ is called a \textit{generalised Cahn-Hoffman field for $u$ at scale $r$} if, for a.e.
$w\in\R^N$, $z_w(C_r^\p)\le1$ and
\beq\label{generalized Cahn Hoffman saturation}
\int_{C_r^\p}\bigl(u(w+x)-u(w+y)\bigr)\,dz_w(x,y) =\osc_{W_r^\p(w)}u.
\eeq
Since the scale $r$ is fixed throughout most of the paper, we will usually omit the expression ``at scale $r$''.
\end{definition}

In the following theorem we consider the case where $g \in \Lip(\R^N)$.

\begin{theorem}\label{existence-uniqueness-Lip}
    Let $h>0$ and let $g \in \Lip(\R^N)$. Then, there exist $z \in L^\infty_{\mathbf{w}}(\R^N ; \mathcal{M}^+(C_r^\p))$ and a function $u_h^g \in \Lip(\R^N)$ such that
    \begin{itemize}
        \item[(a)] for a.e. $w \in \R^N$ one has $z_w(C_r^\p) \leq 1$ and
        \beq\label{Cahn Hoffman equation Lip}
            \int_{C_r^\p} (u_h^g(w+x)-u_h^g(w+y))\,dz_w(x,y)=\mathop{\mathrm{osc}}_{W_r^\p(w)}u_h^g;
        \eeq
        \item[(b)]
        \beq\label{EL equation Lip}
            -h\divr_\p z+u_h^g=g \qquad \text{in }\mathcal M(\R^N).
        \eeq
    \end{itemize}
    Moreover, $\Lip(u_h^g)\leq\Lip(g)$. In addition, if $g\in\mathscr L^+$ (respectively, $g\in\mathscr L^-$), then every Lipschitz solution satisfying (a)-(b) belongs to
    $\mathscr L^+$ (respectively, to $\mathscr L^-$), and the solution is unique in $\Lip(\R^N)$. Finally, if $g$ satisfies $\psi(\nabla g)\leq1$ a.e. for some norm $\psi$, the
    solution may be chosen so that $\psi(\nabla u_h^g)\leq1$ a.e. in $\R^N$.
\end{theorem}

\begin{remark}\label{remark-arbitrary-lip-uniqueness}
The existence part of Theorem~\ref{existence-uniqueness-Lip} requires no coercivity assumption. For general noncoercive Lipschitz data we will not need uniqueness in the whole
class $\Lip(\R^N)$. For coercive data, Step~2 below shows that every Lipschitz solution is automatically in the corresponding class $\mathscr L^+$ or $\mathscr L^-$, so
Proposition~\ref{comparison} yields uniqueness in $\Lip(\R^N)$.
\end{remark}

Before proving Theorem~\ref{existence-uniqueness-Lip}, we establish some properties of generalised Cahn-Hoffman fields that will be used in the proof.

\begin{remark}\label{Remark-Cahn-Hoffman}
Let $(u,z)\in\Lip(\R^N)\times\mathcal M^+(C_r^\p)$, $z(C_r^\p)\le1$, and assume
\beq\label{26}
\int_{C_r^\p}(u(x)-u(y))\,dz(x,y)=\osc_{W_r^\p}u.
\eeq
Then
\beq\label{27}
u(x)=\max_{W_r^\p}u, \qquad u(y)=\min_{W_r^\p}u \quad\text{for }z\text{-a.e. }(x,y)\in C_r^\p.
\eeq
Indeed, the left-hand side of \eqref{26} is bounded above by $z(C_r^\p)\osc_{W_r^\p}u\le\osc_{W_r^\p}u$, so equality implies that the maximum and minimum are attained $z$-a.e. in
the first and second variables, respectively.
\end{remark}

\begin{lemma}\label{char-function-test}
Let $h>0$, $g\in\Lip(\R^N)$, and let $(u,z)$ satisfy (a)-(b) of Theorem~\ref{existence-uniqueness-Lip}. Let $E,F\subset\R^N$ be closed sets such that $\overline{E\triangle F}$ is
compact. Then
\beq\label{char function test equation}
\begin{split}
\frac1{2r}\int_{\R^N}\int_{C_r^\p} &\Big[(\chi_E-\chi_F)(w+x)-(\chi_E-\chi_F)(w+y)\Big]dz_w(x,y)\,dw\\
&=-\frac1h\int_{\R^N}(\chi_E-\chi_F)(u-g)\,dw.
\end{split}
\eeq
\end{lemma}
\begin{proof}
Choose $\Psi\in C_c^\infty(\R^N)$, $0\leq\Psi\leq1$, such that $\Psi\equiv1$ on a neighbourhood of $\overline{E\triangle F}$. For $n\in\N$ set
\[
\vphi_n^E(x):=\bigl(1-n\,\dist(x,E)\bigr)_+, \qquad \vphi_n^F(x):=\bigl(1-n\,\dist(x,F)\bigr)_+,
\]
with the obvious convention if one of the sets is empty, and let $\zeta_n:= (\vphi_n^E-\vphi_n^F) \Psi$. Since $E$ and $F$ are closed, $\vphi_n^E\searrow\chi_E$ and
$\vphi_n^F\searrow\chi_F$ pointwise. Hence $\zeta_n\in\Lip_c(\R^N)$, $|\zeta_n|\leq1$, and $\zeta_n\to\chi_E-\chi_F$ pointwise. Using $\zeta_n$ as a test function in
\eqref{EL equation Lip} gives
\[
\frac1{2r}\int_{\R^N}\int_{C_r^\p} \bigl(\zeta_n(w+x)-\zeta_n(w+y)\bigr)\,dz_w(x,y)\,dw =-\frac1h\int_{\R^N}\zeta_n(u-g)\,dw.
\]
We may pass to the limit by dominated convergence. Indeed, the functions $\zeta_n$ are uniformly bounded and supported in the fixed compact set $\operatorname{supp}\Psi$, while the
integrand on the left vanishes unless $w$ belongs to the fixed compact set $\operatorname{supp}\Psi+W_r^\p$; moreover $z_w(C_r^\p)\leq1$ for a.e. $w$. Letting $n\to\infty$ yields
\eqref{char function test equation}.
\end{proof}

For a Borel set $\Om\subset\R^N$ we define the localised perimeter by
\[
\mathrm{Per}_r^\p(E;\Om) :=\frac1{2r}\int_\Om\osc_{W_r^\p(w)}\chi_E\,dw.
\]

\begin{lemma}\label{minimum-problem-sublevel}
Let $h>0$, $g\in\Lip(\R^N)$, and let $(u,z)$ satisfy (a)-(b) of Theorem~\ref{existence-uniqueness-Lip}. Fix $t\in\R$ and set $E_t:=\{u\leq t\}$. Then, for every $R>r$ and every
closed set $F\subset\R^N$ such that
\[
\overline{E_t\triangle F}\Subset W_{R-r}^\p,
\]
one has
\beq\label{local sublevel minimality}
\begin{split}
\mathrm{Per}_r^\p(E_t;W_R^\p) +\frac1h\int_{E_t\cap W_{R-r}^\p}(g-t)\,dw \leq{}& \mathrm{Per}_r^\p(F;W_R^\p) +\frac1h\int_{F\cap W_{R-r}^\p}(g-t)\,dw.
\end{split}
\eeq
\end{lemma}
\begin{proof}
By Remark~\ref{Remark-Cahn-Hoffman}, for a.e. $w\in\R^N$,
\beq\label{2}
\int_{C_r^\p}\bigl(\chi_{E_t}(w+x)-\chi_{E_t}(w+y)\bigr)\,dz_w(x,y) =-\osc_{W_r^\p(w)}\chi_{E_t}.
\eeq
Indeed, if the oscillation of $\chi_{E_t}$ is one, then $\min_{W_r^\p(w)}u\leq t<\max_{W_r^\p(w)}u$; hence Remark~\ref{Remark-Cahn-Hoffman} implies that $z_w(C_r^\p)=1$ and that
the measure is concentrated on pairs where the first and second variables attain, respectively, the maximum and minimum of $u$. If the oscillation is zero, both sides of \eqref{2}
vanish.

Since $z_w(C_r^\p)\leq1$ for a.e. $w$, we also have
\[
\osc_{W_r^\p(w)}\chi_F \geq -\int_{C_r^\p}\bigl(\chi_F(w+x)-\chi_F(w+y)\bigr)\,dz_w(x,y).
\]
Subtracting \eqref{2} and integrating over $W_R^\p$ yields
\begin{align*}
\mathrm{Per}_r^\p(F;W_R^\p)-\mathrm{Per}_r^\p(E_t;W_R^\p) \geq -\frac1{2r}\int_{W_R^\p}\int_{C_r^\p} &\Big[(\chi_F-\chi_{E_t})(w+x)\\
&-(\chi_F-\chi_{E_t})(w+y)\Big]dz_w(x,y)\,dw.
\end{align*}
The inner integrand vanishes for $w\notin W_R^\p$, because $\overline{E_t\triangle F}\Subset W_{R-r}^\p$ and $x,y\in W_r^\p$. We may therefore extend the outer integral to $\R^N$
and apply Lemma~\ref{char-function-test}. We obtain
\[
\mathrm{Per}_r^\p(F;W_R^\p)-\mathrm{Per}_r^\p(E_t;W_R^\p) \geq \frac1h\int_{\R^N}(\chi_F-\chi_{E_t})(u-g)\,dw.
\]
Since $u>t$ on $F\setminus E_t$ and $u\leq t$ on $E_t\setminus F$, the right-hand side is bounded below by
\[
\frac1h\int_{\R^N}(\chi_F-\chi_{E_t})(t-g)\,dw.
\]
The symmetric difference is contained in $W_{R-r}^\p$, and rearranging the terms gives \eqref{local sublevel minimality}.
\end{proof}

\begin{corollary}\label{bounded-sublevel-minimality-corollary}
Under the assumptions of Lemma~\ref{minimum-problem-sublevel}, if $E_t$ is bounded, then
\beq\label{bounded sublevel minimality}
\mathrm{Per}_r^\p(E_t)+\frac1h\int_{E_t}(g-t)\,dw \leq \mathrm{Per}_r^\p(F)+\frac1h\int_F(g-t)\,dw
\eeq
for every bounded closed $F\subset\R^N$.
\end{corollary}
\begin{proof}
Choose $R>r$ so large that $E_t\cup F\Subset W_{R-r}^\p$. Then the localised perimeters in \eqref{local sublevel minimality} coincide with the corresponding global perimeters, and
the conclusion follows from Lemma~\ref{minimum-problem-sublevel}.
\end{proof}

\begin{proposition}[Comparison principle]\label{comparison}
Let $g_1,g_2\in\Lip(\R^N)$ and, for $i=1,2$, let $u_i\in\Lip(\R^N)$ and $z_i\in L^\infty_{\mathbf w}(\R^N;\mathcal M^+(C_r^\p))$ satisfy {\rm (a)-(b)} of
Theorem~\ref{existence-uniqueness-Lip} with datum $g_i$. Assume that $g_1\ge g_2$. If either \( u_1\in\mathscr L^+\) or \(u_2\in\mathscr L^-\), then \(u_1\ge u_2\).
Moreover, let $g, u\in \Lip(\R^N)$, and $z_u\in L^\infty_{\mathbf w}(\R^N;\mathcal M^+(C_r^\p))$ satisfy {\rm (a)-(b)} of Theorem~\ref{existence-uniqueness-Lip}
with datum $g$. Let $v\in\mathscr L^+$ admit a generalised Cahn-Hoffman field $z_v$ such that
\beq\label{supersolution comparison inequality}
-h\divr_\p z_v+v\ge g \qquad\text{in }\mathcal M(\R^N).
\eeq
Then $u\le v$.
\end{proposition}
\begin{proof}
Assume first that $u_1\in\mathscr L^+$. Fix $t_1<t_2$ and set
\[
E_1:=\{u_1\le t_1\},\qquad E_2:=\{u_2\le t_2\},\qquad D:=E_1\setminus E_2.
\]
Since $E_1$ is bounded and closed, $\overline D$ is compact. Choose $R>r$ so large that $\overline D\Subset W_{R-r}^\p$. Applying Lemma~\ref{minimum-problem-sublevel} to $E_1$ with
competitor $E_1\cap E_2$, and to $E_2$ with competitor $E_1\cup E_2$, gives
\[
\mathrm{Per}_r^\p(E_1;W_R^\p)-\mathrm{Per}_r^\p(E_1\cap E_2;W_R^\p) \le -\frac1h\int_D(g_1-t_1)\,dw,
\]
and
\[
\mathrm{Per}_r^\p(E_2;W_R^\p)-\mathrm{Per}_r^\p(E_1\cup E_2;W_R^\p) \le \frac1h\int_D(g_2-t_2)\,dw.
\]
Adding and using the pointwise submodularity of the oscillation, hence of the localised perimeter, we obtain
\[
0\le \frac1h\int_D\big[(g_2-t_2)-(g_1-t_1)\big] \,dw \le -\frac{t_2-t_1}{h}|D|.
\]
Thus $|D|=0$, and by continuity $E_1\subset E_2$. Since this holds for every $t_1<t_2$, it follows that $u_1\ge u_2$. If $u_2\in\mathscr L^-$, apply the preceding argument to
$-u_2,-u_1$, with the two coordinates in the corresponding Cahn-Hoffman fields exchanged.

For the last assertion we use the monotonicity property of generalised Cahn-Hoffman fields from \cite[Lemma~3.7]{CMR26}; its proof is pointwise and is unchanged when the Euclidean
ball is replaced by $W_r^\p$. In the form needed here, if $a_1,a_2,\Psi:\R\to\R$ are continuous and non-decreasing, then
\[
\int_{C_r^\p}\!\big(q(w+x)-q(w+y)\big)\,d(z_u-z_v)_w(x,y)\ge0, \qquad q:=\Psi\big(a_1(u)-a_2(v)\big),
\]
for a.e. $w$. For $\lambda>0$ choose
\[
q_\lambda:=\big((u\wedge\lambda)-(v\wedge\lambda)\big)^+\in\Lip_c(\R^N).
\]
Subtracting \eqref{supersolution comparison inequality} from the equation for $u$ and testing with $q_\lambda\ge0$ gives
\[
\frac h{2r}\int_{\R^N}\!\int_{C_r^\p} \big(q_\lambda(w+x)-q_\lambda(w+y)\big)\,d(z_u-z_v)_w\,dw +\int_{\R^N}q_\lambda(u-v)\,dw\le0.
\]
The first term is non-negative by the monotonicity property above, with $a_1=a_2=(\,\cdot\,)\wedge\lambda$ and $\Psi=(\,\cdot\,)^+$, while the second one is non-negative because
$q_\lambda>0$ implies $u>v$. Hence $q_\lambda=0$. Letting $\lambda\to+\infty$ yields $u\le v$.
\end{proof}

\begin{proof}[Proof of Theorem~\ref{existence-uniqueness-Lip}]
\textbf{Step 1: existence for arbitrary Lipschitz data.} Let $L:=\Lip(g)$. We approximate $g$ by functions in $\R+L^2$ without increasing $L$. For $R>0$ prescribe $f_R=g$ on
$\overline B_R$ and $f_R=g(0)$ on $\R^N\setminus B_{3R}$; this datum is $L$-Lipschitz on its domain $D:=\overline B_R\cup \big(\R^N\setminus B_{3R}\big)$. Its McShane extension
\beq\label{McShane cutoff}
g_R(x):=\inf_{y\in D}\{f_R(y)+L|x-y|\}
\eeq
satisfies $g_R=g$ on $B_R$, $g_R=g(0)$ outside $B_{3R}$, $\Lip(g_R)\le L$, and $g_R\to g$ locally uniformly as $R\to\infty$. Proposition~\ref{L2-existence-uniqueness} gives
corresponding solutions $(u_R,z^R)$ with $\Lip(u_R)\le L$ and $\|z^R\|_{L^\infty_{\mathbf w}}\le1$.

After passing to a sequence $R_k\to\infty$,
\beq\label{Chan Hoffman weak convergence}
z^{R_k}\overset{*}{\wto}z \quad\text{in }L^\infty_{\mathbf w}(\R^N;\mathcal M(C_r^\p)).
\eeq
The Euler-Lagrange identity tested against any $\vphi\in C_c^\infty$ shows that $u_{R_k}-g_{R_k}$ is locally uniformly bounded: otherwise, by equi-Lipschitz continuity it would
diverge with one sign on a fixed ball, contradicting the boundedness of
\[
\int_{\R^N} \vphi(u_{R_k}-g_{R_k}) \, dw =-\frac h{2r}\int_{\R^N} \!\int_{C^\phi_r}(\vphi(w+x)-\vphi(w+y))\,dz_w^{R_k}(x,y)\,dw.
\]
Thus, up to a further subsequence,
\beq\label{unif converg}
u_{R_k}\to u\qquad\text{locally uniformly},
\eeq
with $\Lip(u)\le L$. Passing to the limit in the Euler-Lagrange identity gives (b). To obtain (a), multiply the Cahn-Hoffman identity for $(u_{R_k},z^{R_k})$ by
$\vphi\in C_c^\infty$ and integrate in $w$. Local uniform convergence gives convergence of the oscillations, while \eqref{Chan Hoffman weak convergence} gives convergence of the
measure term; hence
\[
\int_{\R^N} \vphi(w)\int_{C_r^\p}(u(w+x)-u(w+y))\,dz_w\,dw =\int_{\R^N} \vphi(w)\osc_{W_r^\p(w)}u\,dw,
\]
and (a) follows. This proves existence for every $g\in\Lip(\R^N)$.

\textbf{Step 2: coercive data.} Assume $g\in\mathscr L^+$ and let $u$ be any Lipschitz solution satisfying (a)-(b). We prove that every sublevel set of $u$ is bounded. Set
$L:=\Lip(u)$, fix $t\in\R$, let $s:=t+1$, and choose
\[
\rho:=\frac1{2(1+L)}.
\]
If $x\in E_t:=\{u\leq t\}$, then
\[
u(y)\leq u(x)+L|y-x|<s \qquad\text{for every }y\in\overline B_\rho(x),
\]
so that $\overline B_\rho(x)\subset E_s:=\{u\leq s\}$. Define the closed competitor
\[
F_x:=E_s\setminus B_\rho(x).
\]
Choose $R_x>r$ so large that $\overline B_\rho(x)\Subset W_{R_x-r}^\p$. Applying Lemma~\ref{minimum-problem-sublevel} to $E_s$ and $F_x$ gives
\beq\label{coercivity ball estimate}
\frac1h\int_{B_\rho(x)}(g-s)\,dw \leq \mathrm{Per}_r^\p(F_x;W_{R_x}^\p)-\mathrm{Per}_r^\p(E_s;W_{R_x}^\p).
\eeq
The two characteristic functions coincide outside $B_\rho(x)$. Hence their oscillations on $W_r^\p(w)$ coincide whenever $W_r^\p(w)\cap B_\rho(x)=\emptyset$. Since every
oscillation of a characteristic function belongs to $[0,1]$, the right-hand side of \eqref{coercivity ball estimate} is bounded above by
\[
C_{\rho,r,\p}:=\frac1{2r}\big|B_\rho+W_r^\p\big|,
\]
independently of $x$. Consequently,
\[
|B_\rho|\left(\inf_{B_\rho(x)}g-s\right) \leq hC_{\rho,r,\p} \qquad\text{for every }x\in E_t.
\]
Since $g(x)\to+\infty$ as $|x|\to\infty$, also $\inf_{B_\rho(x)}g\to+\infty$. Thus $E_t$ is bounded. Since $t$ is arbitrary, we deduce that every Lipschitz solution belongs to
$\mathscr L^+$. Existence follows from Step~1, while Proposition~\ref{comparison} gives uniqueness in $\Lip(\R^N)$.

If $g\in\mathscr L^-$, apply the preceding argument to $-g$ and $-u$, exchanging the two coordinates in the Cahn-Hoffman measures. Thus every Lipschitz solution belongs to
$\mathscr L^-$, and uniqueness in $\Lip(\R^N)$ follows in the same way.

\textbf{Step 3: the $\psi$-Lipschitz estimate.} If $\psi(\nabla g)\le1$, repeat the approximation in Step~1 with McShane extensions for the metric $\pso$. Then
$\psi(\nabla g_R)\le1$, and translation invariance and the comparison principle in $L^2$ yield $\psi(\nabla u_R)\le1$. The estimate passes to the locally uniform limit.
\end{proof}

In the following proposition we estimate the solution $u_h^g$ given by Theorem \ref{existence-uniqueness-Lip} when $g$ is of the form $c_2(\po -R)_+ - c_1(\po -R)_-$ for some
$c_2\geq c_1 >0$.

\begin{proposition}\label{Wulff-shape-comparison}
For every $h>0$ and $R>0$, let
\[
 g(w):=c_2(\po(w)-R)_+-c_1(\po(w)-R)_-, \qquad 0<c_1\le c_2.
\]
Let $\eta<0$ be such that
\[
 \beta:=R+\frac{\eta}{c_1}>0.
\]
Then there exist constants $c_0=c_0(r,\eta,R,c_1)>0$ and $h_0=h_0(r,\eta,R,c_1)>0$ such that, for every $0<h\le h_0$,
\beq\label{29}
 u_h^g(w)\le (g(w)+c_0h)\vee\eta \qquad\text{for every }w\in\R^N.
\eeq
The constants can be chosen so that $c_0$ is non-increasing and $h_0$ is non-decreasing as functions of $R$. Moreover, for fixed $\eta,c_1,R$, the constant $c_0$ can be chosen
independently of $r$ when $r>0$ is sufficiently small.
\end{proposition}

\begin{proof}
As in the isotropic case treated in \cite{CMR26}, the proof is based on the construction of a suitable generalised Cahn-Hoffman field. Set
\[
 q(\rho):=c_2(\rho-R)_+-c_1(\rho-R)_-,\qquad v(w):=(q(\po(w))+c_0h)\vee\eta.
\]
Thus $g(w)=q(\po(w))$ and $q$ is non-decreasing. Since $q(\beta)=\eta$, if
\[
 \alpha_h:=\beta-\frac{c_0h}{c_1},
\]
then
\beq\label{plateau v}
 v=\eta\quad\text{on }W_{\alpha_h}^\p, \qquad v-g\ge c_0h\quad\text{in }\R^N.
\eeq
We shall always impose
\beq\label{h small prop49}
 0<h\le \frac{c_1\beta}{8c_0},
\eeq
so that $\alpha_h\ge 7\beta/8$.

It is enough to construct a generalised Cahn-Hoffman field $z$ for $v$ such that, for every non-negative $\zeta\in C_c^0(\R^N)$,
\beq\label{div bound prop49}
 \frac1{2r}\int_{\R^N}\int_{C_r^\p} \bigl(\zeta(w+x)-\zeta(w+y)\bigr)\,dz_w(x,y)\,dw \ge -c_0\int_{\R^N}\zeta\,dw.
\eeq
Indeed, multiplying \eqref{div bound prop49} by $h$ and using \eqref{plateau v} gives
\beq\label{cahn hoffman inequality}
 -h\,\divr_\p z+v\ge g \qquad\text{in }\mathcal M(\R^N).
\eeq
The supersolution comparison in Proposition~\ref{comparison}, applied to \eqref{cahn hoffman inequality}, then yields \eqref{29}.

For $w\ne0$ we write $\widehat w:=w/\po(w)$. At $w=0$ we fix once and for all a vector $e$ with $\po(e)=1$ and put $\widehat 0=e$. Since $v$ is a non-decreasing function of $\po$,
for every $w$ we have
\beq\label{radial extrema prop49}
 \sup_{W_r^\p(w)}v=v(w+r\widehat w).
\eeq
Whenever $\po(w)\ge r$, the corresponding minimum is $v(w-r\widehat w)$.

\medskip
\noindent\textbf{Case 1: $\beta\le 8r$.} Define the following set of \emph{increments}
\beq\label{Dh def prop49}
 D_h(w):=W_r^\p\cap\bigl(W_{\alpha_h}^\p-w\bigr) =\{y\in W_r^\p:\ w+y\in W_{\alpha_h}^\p\}.
\eeq
If $\po(w)<\beta/2+r$, then $D_h(w)$ has positive measure. More precisely,
\beq\label{Dh lower bound}
 |D_h(w)|\ge |W_{\beta/8}^\p|.
\eeq
To see this, observe that the two Wulff balls in \eqref{Dh def prop49} have radii $r$ and $\alpha_h\ge7\beta/8$, while the distance between their centres is less than $r+\beta/2$.
If neither ball contains the other, their intersection contains a Wulff ball of radius
\[
 \frac{r+\alpha_h-\po(w)}2>\frac{\alpha_h-\beta/2}2 \ge\frac{3\beta}{16}>\frac\beta8;
\]
if one contains the other, the smaller radius is at least $\beta/8$ because $\beta\le8r$. This proves \eqref{Dh lower bound}.

We now set
\beq\label{z general prop49}
 z_w:=
 \begin{cases}
 \de_{r\widehat w}\times\de_{-r\widehat w}, & \po(w)\ge\beta/2+r,\\[1mm]
 \displaystyle \de_{r\widehat w}\times \frac{\mathcal L^N\res D_h(w)}{|D_h(w)|}, & \po(w)<\beta/2+r.
 \end{cases}
\eeq
The family $w\mapsto z_w$ is weakly measurable and $z_w(C_r^\p)=1$. We next verify the Cahn-Hoffman identity. In the first region $\po(w)\ge\beta/2+r>r$, so by radial monotonicity
\[
 \osc_{W_r^\p(w)}v =v(w+r\widehat w)-v(w-r\widehat w),
\]
which is exactly the integral against the first line of \eqref{z general prop49}. In the second region, every $y\in D_h(w)$ satisfies $v(w+y)=\eta$ by \eqref{plateau v}. Since
$v\ge\eta$ everywhere and \eqref{radial extrema prop49} gives the maximum, we obtain
\[
 \osc_{W_r^\p(w)}v =v(w+r\widehat w)-\eta =\int_{C_r^\p}(v(w+x)-v(w+y))\,dz_w(x,y).
\]
Thus $z$ is a generalised Cahn-Hoffman field for $v$.

It remains to estimate the divergence. Split the left-hand side of \eqref{div bound prop49} according to the two regions in \eqref{z general prop49}.

For the first region, since $\zeta\geq0$, we may discard the non-negative term and obtain
\[
\begin{aligned}
&\frac1{2r}\int_{\{\po(w)\geq\beta/2+r\}} \bigl(\zeta(w+r\widehat w)-\zeta(w-r\widehat w)\bigr)\,dw\\
&\qquad\geq -\frac1{2r}\int_{\{\po(w)\geq\beta/2+r\}} \zeta(w-r\widehat w)\,dw\\
&\qquad= -\frac1{2r}\int_{\beta/2}^{+\infty} \left(1+\frac rs\right)^{N-1} \int_{\{\po(w)=s\}} \zeta(w)\frac{d\H^{N-1}(w)} {|\nabla\po(w)|}\,ds.
\end{aligned}
\]
Here in the last equality we used the coarea formula and the change of variables $s=\rho-r$, together with the $0$-homogeneity of $|\nabla\po|$. Since $s\geq\beta/2$,
\[
\left(1+\frac rs\right)^{N-1} \leq \left(1+\frac{2r}{\beta}\right)^{N-1},
\]
and therefore, by the coarea formula,
\[
\begin{aligned}
&\frac1{2r}\int_{\{\po(w)\geq\beta/2+r\}} \bigl(\zeta(w+r\widehat w)-\zeta(w-r\widehat w)\bigr)\,dw\\
&\qquad\geq -\frac1{2r} \left(1+\frac{2r}{\beta}\right)^{N-1} \int_{\R^N}\zeta(w)\,dw.
\end{aligned}
\]

For the second region, we discard the non-negative term containing $\zeta(w+r\widehat w)$. Using \eqref{Dh lower bound}, $D_h(w)\subset W_r^\p$, and $\zeta\geq0$, we obtain
\[
\begin{aligned}
&-\frac1{2r}\int_{\{\po(w)<\beta/2+r\}} \frac1{|D_h(w)|}\int_{D_h(w)}\zeta(w+y)\,dy\,dw\\
&\qquad\geq -\frac1{2r}\frac1{|W_{\beta/8}^\p|} \int_{\R^N}\int_{W_r^\p}\zeta(w+y)\,dy\,dw\\
&\qquad= -\frac1{2r}\frac{|W_r^\p|}{|W_{\beta/8}^\p|} \int_{\R^N}\zeta(w)\,dw\\
&\qquad= -\frac1{2r}\left(\frac{8r}{\beta}\right)^N \int_{\R^N}\zeta(w)\,dw.
\end{aligned}
\]
Hence \eqref{div bound prop49} holds in this case with $c_0$ replaced by
\beq\label{cgen prop49}
c(r,\beta):= \frac1{2r}\left[ \left(1+\frac{2r}{\beta}\right)^{N-1} +\left(\frac{8r}{\beta}\right)^N \right].
\eeq
\medskip
\noindent\textbf{Case 2: $\beta\ge 8r$.} Let
\[
 t(\rho):=
 \begin{cases}
 0,&0\le\rho\le\beta/2,\\
 4\rho/\beta-2,&\beta/2<\rho<3\beta/4,\\
 1,&\rho\ge3\beta/4,
 \end{cases}
\]
and set
\beq\label{z refined prop49}
 z_w:=\de_{r t(\po(w))\widehat w} \times\de_{-r t(\po(w))\widehat w}.
\eeq
If $\po(w)\le3\beta/4$, then, using $r\le\beta/8$ and \eqref{h small prop49}, every $x\in W_r^\p(w)$ satisfies $\po(x)\le7\beta/8\le\alpha_h$; hence $v\equiv\eta$ on $W_r^\p(w)$
and the Cahn-Hoffman identity is trivial. If $\po(w)>3\beta/4$, then $t(\po(w))=1$ and the Cahn-Hoffman identity follows from the radial extrema exactly as above. Thus
\eqref{z refined prop49} is again a generalised Cahn-Hoffman field for $v$.

For $\zeta\ge0$ set
\[
 A_\zeta(s):=\int_{\{\po=s\}} \zeta(w)\,\frac{d\H^{N-1}(w)}{|\nabla\po(w)|}.
\]
Arguing as above, using the $0$-homogeneity of $|\nabla\po|$, the coarea formula and a change of variable, after a lengthy but straightforward calculation, we obtain
\[
 \int_{\R^N}\zeta\bigl(w\pm rt(\po(w))\widehat w\bigr)\,dw =\int_0^\infty a_\pm(s)A_\zeta(s)\,ds,
\]
where
\[
 a_+(s)=
 \begin{cases}
 1,&s\le\beta/2,\\[1mm]
 \displaystyle\frac{(1+2r/s)^{N-1}}{(1+4r/\beta)^N}, &\beta/2<s<3\beta/4+r,\\[2mm]
 (1-r/s)^{N-1},&s\ge3\beta/4+r,
 \end{cases}
\]
and
\[
 a_-(s)=
 \begin{cases}
 1,&s\le\beta/2,\\[1mm]
 \displaystyle\frac{(1-2r/s)^{N-1}}{(1-4r/\beta)^N}, &\beta/2<s<3\beta/4-r,\\[2mm]
 (1+r/s)^{N-1},&s\ge3\beta/4-r.
 \end{cases}
\]
 We now estimate the difference $a_+-a_-$. Set \( \e:=\frac r\beta,\, x:=\frac rs. \) Since $\beta\ge 8r$ and $s\ge\beta/2$, we have \( 0\le x\le2\e,\, 4\e\le\frac12. \) Let
\[
f(t):=(1+t)^{-N},\qquad g(t):=(1+t)^{N-1},
\]
and let $L_N$ be a common Lipschitz constant for $f$ and $g$ on $[-1/2,1/2]$. We claim that, for a constant $C_N>0$ depending only on $N$,
\beq\label{density bound prop49}
a_+(s)-a_-(s)\ge -C_N\frac r\beta \qquad\text{for every }s\ge\frac\beta2.
\eeq
Indeed, consider for instance the range \( \frac\beta2<s<\frac{3\beta}{4}-r. \) In this case, \(a_+(s)-a_-(s) =f(4\e)g(2x)-f(-4\e)g(-2x), \) and hence
\[
a_+(s)-a_-(s) = f(4\e)\bigl(g(2x)-g(-2x)\bigr)\\
+ \bigl(f(4\e)-f(-4\e)\bigr)g(-2x).
\]
The first term is non-negative, since $g$ is non-decreasing, while $0\le g(-2x)\le1$. Therefore
\[
a_+(s)-a_-(s)\ge -8L_N\e =-8L_N\frac r\beta.
\]
The other two ranges determined by $3\beta/4-r$ and $3\beta/4+r$ are treated analogously (and more directly), using the Lipschitz continuity of $f$ and $g$ and the monotonicity of
$g$. This proves \eqref{density bound prop49}, after possibly enlarging $C_N$. Consequently, since $A_\zeta\ge0$,
\begin{multline*}
\frac1{2r}\int_{\R^N}\int_{C_r^\p} \bigl(\zeta(w+x)-\zeta(w+y)\bigr)\,dz_w(x,y)\,dw\\
= \frac1{2r}\int_{\beta/2}^{+\infty} (a_+(s)-a_-(s))A_\zeta(s)\,ds\geq -\frac{C_N}{2\beta} \int_{\beta/2}^{+\infty}A_\zeta(s)\,ds \ge -\frac{C_N}{2\beta}\int_{\R^N}\zeta\,dw=:
-\frac{K_N}{\beta}\int_{\R^N}\zeta\,dw
\end{multline*}
In particular, this bound is independent of $r$.

\medskip
We finally choose the constants so that they have the required monotonicity properties. Set
\[
 c_0=c_0(r,\beta):=
 \begin{cases}
 \max\{c(r,\beta),K_N/\beta\},&\beta\le8r,\\[1mm]
 K_N/\beta,&\beta\ge8r,
 \end{cases}
 \qquad h_0(r,\beta,c_1):=\frac{c_1\beta}{8c_0(r,\beta)}.
\]
Note that enlarging $K_N$ if necessary, we may also assume that
\[
c(r,8r)\leq \frac{K_N}{8r}.
\]
This guarantees that the two definitions above match at $\beta=8r$ and that $c_0$ is non-increasing in $\beta$ (and therefore in $R$), while $h_0$ is non-decreasing. The choice of
$h_0$ gives exactly \eqref{h small prop49}, so all the previous constructions apply. Finally, for fixed $\eta,c_1,R$ we have $\beta$ fixed, and if $r\le\beta/8$ then
$c_0=K_N/\beta$, which is independent of $r$. This concludes the proof.
\end{proof}

\begin{proposition}\label{comparison-phi-Wulff-shape}
    For every $h>0$, let $\theta^h = u_h^g$ when $g = \po$. Then, for every $\eta >0$ there exist $c_0 =c_0(r,\eta) >0$ and $h_0 = h_0(r,\eta) >0$ such that for every $h \leq h_0$
    \beq
        \theta^h (w) \leq (\po (w) + c_0h) \vee \eta \hspace{5mm} \text{for every } w \in \R^N.
    \eeq
    Moreover, for every $\alpha \in (0, 1 /(N+1))$ there exist $h_1=h_1(\alpha)$ and a constant $c>0$ such that
    \beq\label{estimate zero Wulff shape}
        \theta^h(0) \leq c \, h^\alpha
    \eeq
    for every $h \leq h_1$.
\end{proposition}
\begin{proof}
    The first part follows by arguing as in the proof of Proposition \ref{Wulff-shape-comparison}, with $\eta/c_1 + R$ replaced by $\eta$. As for \eqref{estimate zero Wulff shape},
    the proof is the same as in \cite[Prop.~3.11]{CMR26}.
\end{proof}

\subsection{Weak formulation of the flow.} In this section we introduce the weak formulation of the anisotropic Minkowski flow \eqref{curvature}, with $\mathcal{K}_r^\p$ denoting
the anisotropic nonlocal curvature associated with $\mathrm{Per}_r^\p$. The definition below extends \cite[Definition~3.12]{CMR26} to the anisotropic setting and to initial sets
whose boundary need not be compact. As in \cite{CMP17,CMNP19}, it is formulated through a distributional inequality for the distance function. Throughout the paper the norms $\p$
and $\psi$ will stand for the \textit{anisotropy} and the \textit{mobility}, respectively. We remark that we do not assume any regularity on $\p$ (nor on $\psi$). We observe that
there exist two positive constants $c_1 \leq c_2$ such that
\[
c_1 \po \leq \pso \leq c_2 \po.
\]

Given $z \in L_{\mathbf{w}}^\infty( \R^N \times (0,T) ; \mathcal{M}^+(C_r^\p))$ we still denote by $\divr_\p z$ its anisotropic nonlocal divergence with respect to $\p$, that is
the measure on $\R^N \times (0,T)$ satisfying
\beq
    \int_{\R^N \times (0,T)} \vphi (w,t) \, d ( \divr_\p z)(w,t) = - \frac{1}{2r} \int_0^T \int_{\R^N} \int_{C_r^\p} (\vphi (w+x,t) - \vphi (w+y,t)) \, dz_{(w,t)}(x,y) dw dt
\eeq
for every $\vphi \in C_c^0(\R^N \times (0,T))$, where $z_{(w,t)}$ denotes the evaluation at $(w,t)$ of the map $z$. Note that the above integrals are also well defined for
$\vphi \in L^1(0,T;C^0(K))$ for any compact set $K \subset \R^N$.

We recall that, given a set $E$ and $\de >0$, we denote by $(E)_\de$ its anisotropic $\de$-neighbourhood with respect to $\p$, namely the set $E + W_\de^\p$. We also recall that a
sequence of closed sets $(E_n)_{n \in \N}$ of $\R^d$ converges to a closed set $E \subset \R^d$ in the Kuratowski sense, and we write $E_n \overset{\mathscr{K}}{\rightarrow} E$, if
the following two conditions hold:
\begin{itemize}
    \item[(i)] if $x_n \in E_n$ for all $n$, then any cluster point of $( x_n )_{n \in \N}$ belongs to $E$;
    \item[(ii)] any $x \in E$ is the limit of a sequence $(x_n)_{n \in \N}$, with $x_n \in E_n$ for all $n$.
\end{itemize}
In the following, given a subset $E$ of $\R^N \times [0, + \infty)$, with $\mathrm{Int}(E) $ we denote the interior of $E$ in the relative topology of $\R^N \times [0, + \infty)$.
\begin{definition}\label{super-flow-def}
    Let $E^0 \subset \R^N$ be a closed set. Let $E$ be a closed set in $\R^N \times [0, + \infty)$ and for each $t \geq 0$ denote $E(t) \coloneqq \{ x \in \R^N : (x,t) \in E \}$.
    We say that $E$ is a $\textit{weak superflow }$ of the curvature flow \eqref{curvature} with initial datum $E^0$ if
    \begin{itemize}
        \item[(a)] $E(0) \subset E^0$;
        \item[(b)] for every $\de > 0$ there exists $\tau(\de) > 0$ with the following property: for every $t \geq 0$ and for every $s \in (t,t+ \tau(\de))$ we have $E(s) \subset (E(t))_\de$; 
        \item[(c)] $E(s) \overset{\mathscr{K}}{\rightarrow} E(t)$ as $s \nearrow t$ for all $t >0$;
        \item[(d)] setting $T^{\mathrm{ext}} \coloneqq \inf \{ t > 0 \, : \, E(s) = \emptyset \text{ for } s \geq t\}$ and $d(x,t) \coloneqq \dist^{\pso} (x,E(t))$ for every
        $(x,t) \in \R^N \times (0, T^{\mathrm{ext}})$, there exists $z \in L^\infty_{\mathbf{w}}(\R^N \times (0, T^{\mathrm{ext}}); \mathcal{M}^+(C_r^\p))$ such that the following
        holds:
        \begin{itemize}
            \item[(d1)] we have
\beq\label{24}
    \pa_t d \geq \divr_\p z \hspace{5mm} \text{as measures in } (\R^N \times (0, T^{\mathrm{ext}})) \setminus E;
\eeq
\item[(d2)] for a.e. $(w,t) \in \R^N \times (0, T^{\mathrm{ext}})$, we have $z_{(w,t)}(C_r^\p) \leq 1$ and
\beq\label{CH}
    \int_{C_r^\p} (d(w+x,t) - d(w+y,t)) \, dz_{(w,t)}(x,y) = \mathop{\mathrm{osc}}_{W_r^\p(w)} d( \cdot, t);
\eeq
\item[(d3)] for every $\lambda > 0$
\beq\label{25}
    ( \divr_\p z)^+ \in L^\infty ( \{ (w,t) \in \R^N \times (0, T^{\mathrm{ext}}) \, : \, d(w,t) \geq \lambda \}).
\eeq
         \end{itemize}
    \end{itemize}
Let $A$ be a relatively open set in $\R^N \times [0, + \infty)$. We say that $A$ is a $\textit{weak subflow }$ of the curvature flow \eqref{curvature}, with initial datum $E^0$ if
$(\R^N \times [0, + \infty)) \setminus A$ is a superflow with initial datum $\R^N \setminus \mathrm{Int}_{\R^N} (E^0)$, where $\mathrm{Int}_{\R^N}(E^0)$ denotes the interior of
$E^0$ in the topology of $\R^N$.

Finally, let $E^0 \subset \R^N$ be a closed set, coinciding with the closure of its interior. Then, we say that a closed set $E \subset \R^N \times [0, + \infty)$ is a
$\textit{weak flow }$ with initial datum $E^0$ if $E = \overline{\mathrm{Int}(E)}$, the set $E$ is a weak superflow with initial datum $E^0$, and $\mathrm{Int}(E)$ is a weak
subflow with initial datum $E^0$.
\end{definition}
\begin{remark}\label{exitiction-time-remark}
    Note that, in particular, taking into account that $(\emptyset)_\de = \emptyset$ for every $\de >0$, condition (b) above implies the following:
    \begin{itemize}
        \item[(b')] if $t \geq 0$ and $E(t) = \emptyset$, then $E(s) = \emptyset$ for all $s>t$.
    \end{itemize}
\end{remark}
We now define the weak level set flows associated with \eqref{curvature}, following \cite[Definition~3.14]{CMR26} and the corresponding local formulation in \cite{CMNP19}.

\begin{definition}[Level set subsolutions and supersolutions]\label{sup-sub-levle-set-def} Let $u^0$ be a uniformly continuous function on $\R^N$ such that $\{u^0\leq\lambda\}$ has compact boundary for every $\lambda\in\R$. We will say that a lower semicontinuous function $u \colon \R^N \times [0, + \infty) \rightarrow \R$ is a $\textit{level set supersolution}$ corresponding to \eqref{curvature}, with initial datum $u^0$, if $u(\cdot , 0) \geq u^0( \cdot)$ and if for a.e. $\lambda \in \R$ the closed sublevel set $\{(x,t) \, : \, u(x,t) \leq \lambda \}$ is a superflow of \eqref{curvature} in the sense of Definition \ref{super-flow-def}, with initial datum $\{u^0 \leq \lambda \}$.

We will say that an upper-semicontinuous function $u \colon \R^N \times [0, + \infty) \rightarrow \R$ is a $\textit{level set subsolution}$ corresponding to \eqref{curvature}, with
initial datum $u^0$, if $-u$ is a level set supersolution in the previous sense, with initial datum $-u^0$.

Finally, we will say that a continuous function $u \colon \R^N \times [0, + \infty) \rightarrow \R$ is a $\textit{level set solution}$ to \eqref{curvature} if it is both a
subsolution and a supersolution.
\end{definition}

\subsection{Minimising movements.}\label{min-mov-scheme-section} In this section we provide an existence result for the weak flows introduced in Definition \ref{super-flow-def}, by
adapting to the nonlocal setting a variant of the Almgren-Taylor-Wang minimising movement scheme (see \cite{ATW}) introduced in \cite{CMP12} and \cite{Chambolle}.

\vspace{2mm}

Let $E^0 \subset \R^N$ be closed. Fix a time step $h>0$ and set $E_h^0:=E^0$. We inductively define $E_h^{k+1}$ as follows. If $E_h^k\neq\emptyset,\R^N$, the signed distance
$d_{E_h^k}^{\pso}$ is $\pso$-Lipschitz, and Theorem~\ref{existence-uniqueness-Lip} provides at least one pair
\[
(u_h^{k+1},z_h^{k+1})\in\Lip(\R^N)\times L_{\mathbf w}^\infty(\R^N;\mathcal M^+(C_r^\p))
\]
satisfying
\beq\label{min mov scheme}
    \begin{cases}
        \displaystyle \int_{C_r^\p}(u_h^{k+1}(w+x) - u_h^{k+1}(w+y)) \, d(z_h^{k+1})_w(x,y) = \mathop{\mathrm{osc}}_{W_r^\p(w)} u_h^{k+1} & \text{for a.e. } w \in \R^N; \\[2pt]
        (z_h^{k+1})_w (C_r^\p) \leq 1, & \text{for a.e. } w \in \R^N; \\[2pt]
        - h \divr_\p z_h^{k+1} + u_h^{k+1} = d^{\pso}_{E_h^k} & \text{in } \mathcal{M}(\R^N),
    \end{cases}
\eeq
and such that $\psi(\nabla u_h^{k+1})\le1$ a.e. We choose one such pair and set $E_h^{k+1}:=\{u_h^{k+1}\le0\}$. If either $E_h^k=\emptyset$ or $E_h^k=\R^N$, we set
$E_h^{k+1}:=E_h^k$. Thus the construction is well defined for every closed initial set, although uniqueness is not guaranteed when the signed-distance datum is
noncoercive. If $E_h^k$ or its complement is bounded, the datum belongs to $\mathscr L^+$ or $\mathscr L^-$, respectively, and Theorem~\ref{existence-uniqueness-Lip} gives
uniqueness at that step.

We denote by $T_h^{\mathrm{ext}}$ the \textit{discrete extinction time}, that is the first discrete time $hk$ such that $E_h^k=\emptyset$, if such a time exists; otherwise
$T_h^{\mathrm{ext}}=+\infty$. Analogously, let $T_h^{\mathrm{all}}$ be the first discrete time $hk$ such that $E_h^k=\R^N$, if such a time exists, and $+\infty$ otherwise. Finally,
set
\[
T_h^*:=\min\{T_h^{\mathrm{ext}},T_h^{\mathrm{all}}\}.
\]

Since $\psi(\nabla u_h^{k+1})\le1$, for every $\tau\in\R^N$,
\[
u_h^{k+1}(\cdot+\tau)-u_h^{k+1}\le\pso(\tau).
\]
Hence, if $x\in E_h^{k+1}$,
\[
u_h^{k+1}(w)\le u_h^{k+1}(w)-u_h^{k+1}(x)\le\pso(w-x),
\]
and minimising over $x\in E_h^{k+1}$ gives
\beq\label{min mov scheme inequality1}
    u_h^{k+1}(w) \leq d^{\pso}_{E_h^{k+1}}(w) \qquad \text{for every } w \in \{d^{\pso}_{E_h^{k+1}}>0\}.
\eeq
Similarly, if $w\in\{d^{\pso}_{E_h^{k+1}}<0\}$ and $x\notin E_h^{k+1}$, then $u_h^{k+1}(x)>0$ and
\[
u_h^{k+1}(w)\ge u_h^{k+1}(x)-\pso(w-x)>-\pso(w-x).
\]
Taking the supremum over $x\notin E_h^{k+1}$ yields
\beq\label{min mov scheme inequality2}
    u_h^{k+1}(w) \geq d^{\pso}_{E_h^{k+1}}(w) \qquad \text{for every } w \in \{d^{\pso}_{E_h^{k+1}}<0\}.
\eeq

We are now in a position to define the time discrete evolutions. For every $h>0$ we set
\beq\label{discrete evolution def}
\begin{split}
    & E_h := \{ (w,t) \in \R^N \times [0,+\infty) : w \in E_h^{[t/h]} \}, \\
    & E_h(t) := E_h^{[t/h]} = \{ w \in \R^N : (w,t) \in E_h \}, \qquad t \geq 0, \\
    & d_h(\cdot,t) := d_{E_h(t)}^{\pso}, \qquad d_h : \R^N \times [0,+\infty) \rightarrow \overline{\R}, \\
    & u_h(\cdot,t) := u_h^{[t/h]}, \qquad u_h : \R^N \times [h,T_h^*) \rightarrow \R, \\
    & z_h(\cdot,t) := z_h^{[t/h]}, \qquad z_h : \R^N \times [h,T_h^*) \rightarrow \mathcal{M}^+(C_r^\p),
    \end{split}
\eeq
where $[\cdot]$ stands for the integer part of its argument.

\begin{remark}\label{Discr-avoid-remark}[Discrete comparison, avoidance, and Wulff barriers]
Assume first that $E^0,F^0\subset\R^N$ are closed sets such that both $d_{E^0}^{\pso}$ and $d_{F^0}^{\pso}$ are coercive, and let $E_h$ and $F_h$ be the corresponding discrete
evolutions. Then the scheme satisfies the discrete comparison principle: if $E^0\subset F^0$, then
\[
E_h(t)\subset F_h(t)\qquad\text{for every }t\ge0.
\]
Indeed, $E^0\subset F^0$ implies $d_{E^0}^{\pso}\ge d_{F^0}^{\pso}$, and Proposition~\ref{comparison} yields the corresponding ordering of the solutions of the discrete scheme. The
conclusion follows by iteration. Notice that coercivity is preserved at every step, so that the discrete evolutions are uniquely determined.

The scheme also satisfies the discrete avoidance principle. Namely, if $\dist^{\pso}(E^0,F^0)\ge\Delta>0$, then
\[
\dist^{\pso}(E_h(t),F_h(t))\ge\Delta \qquad\text{for every }t\ge0.
\]
Indeed, the initial separation is equivalent to $d_{E^0}^{\pso}+d_{F^0}^{\pso}\ge\Delta$. Hence, by Proposition~\ref{comparison}, together with sign symmetry and invariance under
addition of constants, the corresponding solutions at the first step satisfy $u_E^1\ge\Delta-u_F^1$. If $x\in E_h^1$, then $u_E^1(x)\le0$ and therefore $u_F^1(x)\ge\Delta$. By
\eqref{min mov scheme inequality1}, $d_{F_h^1}^{\pso}(x)\ge\Delta$. Thus $\dist^{\pso}(E_h^1,F_h^1)\ge\Delta$, and the conclusion follows by iteration.

For a general closed initial set the discrete evolution need not be unique, and the preceding comparison principle is not available in full generality. Nevertheless, comparison
with bounded Wulff barriers still holds. More precisely, let $E_h$ be any discrete evolution starting from a closed set $E^0$, let $W^0=W_R^\psi(x_0)$, and let $W_h$ be the unique
discrete evolution starting from $W^0$. If $W^0\subset E^0$, then
\beq\label{one sided discrete inclusion}
W_h(t)\subset E_h(t)\qquad\text{for every }t\ge0.
\eeq
Indeed, $d_{W^0}^{\pso}\ge d_{E^0}^{\pso}$ and the solution corresponding to $d_{W^0}^{\pso}$ belongs to $\mathscr L^+$, so the one-sided comparison in Proposition~\ref{comparison}
applies; the argument can then be iterated.

Similarly, if $\dist^{\pso}(E^0,W^0)\ge\Delta>0$, then
\beq\label{one sided discrete avoidance}
\dist^{\pso}(E_h(t),W_h(t))\ge\Delta \qquad\text{for every }t\ge0.
\eeq
Indeed, $d_{E^0}^{\pso}\ge\Delta-d_{W^0}^{\pso}$ and $\Delta-d_{W^0}^{\pso}\in\mathscr L^-$. Hence the one-sided comparison in Proposition~\ref{comparison} gives
$u_E^1\ge\Delta-u_W^1$. Arguing as above and iterating proves \eqref{one sided discrete avoidance}.
\end{remark}
\subsection{Evolution of $\p$-Wulff shapes}
In this section we derive an inner barrier for the minimising movements starting from a $\p$-Wulff shape. We recall that $0<c_1\le c_2$ are chosen so that
\beq\label{costants psi phi}
    c_1\po\le\pso\le c_2\po.
\eeq
In particular, for every $\rho>0$,
\beq\label{44}
    d^{\pso}_{W^\p_\rho} \le c_2(\po-\rho)_+-c_1(\po-\rho)_-.
\eeq

\begin{proposition}\label{evolution-phi-wulff-shape}
For every $R>0$ and every $0<\de<R$ there exist $\gamma_0=\gamma_0(r,\de,c_1)>0$ and $h_0=h_0(r,\de,c_1)>0$ such that, if $E^0=W_R^\p$, then
\[
E_h(s)\supseteq W^\p_{R-\gamma_0s}
\]
for every $0<h\le h_0$ and every $0<s<(R-\de)/\gamma_0$. Moreover, for $r>0$ sufficiently small, $\gamma_0$ and $h_0$ can be chosen independently of $r$.
\end{proposition}

\begin{proof}
Fix $0<\de<R$ and set $\eta:=-c_1\de/2$. For $\rho\ge\de$ let
\[
q_\rho(w):=c_2(\po(w)-\rho)_+-c_1(\po(w)-\rho)_-.
\]
Apply Proposition~\ref{Wulff-shape-comparison} with this fixed value of $\eta$. Denote by $c_*$ and $h_*$ the constants provided there for the reference radius $\rho=\de$. Since,
for fixed $\eta$, the constant $c_0$ is non-increasing and $h_0$ is non-decreasing with respect to the radius, for every $\rho\ge\de$ and $0<h\le h_*$ the solution
$\widetilde v_h^\rho$ corresponding to the datum $q_\rho$ satisfies
\beq\label{phi wulff one step upper}
\widetilde v_h^\rho(w) \le \bigl(q_\rho(w)+c_*h\bigr)\vee\eta \qquad\text{for every }w\in\R^N.
\eeq
After reducing $h_*$, if necessary, we may also assume $c_*h_*/c_1<\de/2$.

Let $u_h^\rho$ be the solution corresponding to the datum $d^{\pso}_{W_\rho^\p}$. By \eqref{44} and Proposition~\ref{comparison}, $u_h^\rho\le\widetilde v_h^\rho$. Hence
\eqref{phi wulff one step upper} gives
\beq\label{phi wulff one step}
\{u_h^\rho\le0\}\supseteq W^\p_{\rho-\gamma_0h}, \qquad \gamma_0:=\frac{c_*}{c_1},
\eeq
for every $\rho\ge\de$ and $0<h\le h_*$. Indeed, on $W^\p_{\rho-c_*h/c_1}$ one has $q_\rho+c_*h\le0$, while $\eta<0$.

We now iterate \eqref{phi wulff one step}. If $\rho_k:=R-\gamma_0kh>\de$ and $E_h^k\supseteq W^\p_{\rho_k}$, the discrete comparison principle implies
$E_h^{k+1}\supseteq W^\p_{\rho_{k+1}}$. Therefore, for $0<s<(R-\de)/\gamma_0$,
\[
E_h(s)=E_h^{[s/h]} \supseteq W^\p_{R-\gamma_0[s/h]h} \supseteq W^\p_{R-\gamma_0s}.
\]
This proves the claim with $h_0:=h_*$. Finally, for the fixed reference radius $\de$ and the fixed value $\eta=-c_1\de/2$, the last part of Proposition~\ref{Wulff-shape-comparison}
(together with the explicit choice of $h_0$ made in its proof) allows $c_*$, and hence both $\gamma_0$ and $h_0$, to be chosen independently of $r$ when $r$ is sufficiently small.
\end{proof}

\subsection{Evolution of $\psi$-Wulff shapes}
From this subsection on, whenever $\psi$-Wulff barriers are used, we assume that $\psi$ is $\p$-regular in the sense below.

\begin{definition}\label{psi-phi regularity}
Following \cite{CMNP19}, we say that a norm $\psi$ is $\p$-regular if there exists $\e_0 \in (0,1)$ such that, for every $w\in\pa W_1^\psi$, there is $y\in W_1^\psi$ satisfying
\[
W_{\e_0}^\p(y)\subset W_1^\psi, \qquad w\in\pa W_{\e_0}^\p(y).
\]
\end{definition}

Since $\p$ and $\psi$ are fixed, in what follows we omit the dependence of the constants on $c_1,c_2$ and on the $\p$-regularity constant $\e_0$.

\begin{proposition}\label{psi-Wulff-shape-evolution}
Assume that $\psi$ is $\p$-regular. For every $R>0$ and every $0<\de<R$ there exist $\gamma_0=\gamma_0(r,\de)>0$ and $h_0=h_0(r,\de)>0$ such that, if $E^0=W_R^\psi$, then
\[
E_h(s)\supseteq W^\psi_{R-\gamma_0s}
\]
for every $0<h\le h_0$ and every $0<s<(R-\de)/\gamma_0$. Moreover, for $r>0$ sufficiently small, $\gamma_0$ can be chosen independently of $r$.
\end{proposition}

\begin{proof}
Fix $0<\de<R$ and set $\de_\p:=\e_0\de/2$. Let $\gamma_\p=\gamma_\p(r,\de_\p,c_1)$ and $h_\p=h_\p(r,\de_\p,c_1)$ be the constants given by
Proposition~\ref{evolution-phi-wulff-shape}. Reducing $h_\p$ if necessary, assume also $\gamma_\p h_\p<\de_\p$ and $c_2\gamma_\p h_\p<\de/2$.

For $\rho\ge\de$ set  
\[
C_\rho:=\{y:\ W_{\e_0\rho}^\p(y)\subset W_\rho^\psi\}.
\]
Scaling Definition~\ref{psi-phi regularity}, every point of $\pa W_\rho^\psi$ is touched from the inside by a translate of $W_{\e_0\rho}^\p$. Hence
\[
W_\rho^\psi=C_\rho+W_{\e_0\rho}^\p,
\]
since the right-hand side is a convex subset of $W_\rho^\psi$ containing its boundary. Let $0\le a<\e_0\rho$. We claim that
\beq\label{rolling wulff erosion}
W^\psi_{\rho-c_2a} \subset C_\rho+W^\p_{\e_0\rho-a}.
\eeq
Indeed, let $x\in W^\psi_{\rho-c_2a}$. If $v\in W_a^\p$, then $\pso(x+v)\le\pso(x)+\pso(v) \le\rho-c_2a+c_2\po(v)\le\rho$. Hence $x+W_a^\p\subset W_\rho^\psi=C_\rho+W_{\e_0\rho}^\p$.
Fix $\xi\in\R^N$ and take the supremum of the scalar product with $\xi$ over the two sets in the previous inclusion. Since $\sup_{v\in W_s^\p}v\cdot\xi=s\p(\xi)$ and, for a
Minkowski sum, the two variables vary independently, we obtain
\[
x\cdot\xi+a\p(\xi) \le \sup_{y\in C_\rho}y\cdot\xi+\e_0\rho\p(\xi), \qquad\text{hence}\qquad x\cdot\xi \le \sup_{z\in C_\rho+W_{\e_0\rho-a}^\p}z\cdot\xi.
\]
The latter inequality holds for every $\xi\in\R^N$. Since $C_\rho+W_{\e_0\rho-a}^\p$ is closed and convex, we deduce $x\in C_\rho+W_{\e_0\rho-a}^\p$. This proves
\eqref{rolling wulff erosion}.

Consider one step of the scheme starting from $W_\rho^\psi$, with $\rho\ge\de$. For every $y$ such that $W_{\e_0\rho}^\p(y)\subset W_\rho^\psi$,
Proposition~\ref{evolution-phi-wulff-shape}, applied at time $h$, and the discrete comparison principle give
\[
E_h^1\supseteq W_{\e_0\rho-\gamma_\p h}^\p(y) \qquad(0<h\le h_\p).
\]
Taking the union in $y$ and using \eqref{rolling wulff erosion} with $a=\gamma_\p h$ yields
\beq\label{psi wulff one step}
E_h^1\supseteq W^\psi_{\rho-c_2\gamma_\p h}.
\eeq
Thus, with $\gamma_0:=c_2\gamma_\p$, the same induction as in the proof of Proposition~\ref{evolution-phi-wulff-shape} gives
\[
E_h(s)\supseteq W^\psi_{R-\gamma_0[s/h]h} \supseteq W^\psi_{R-\gamma_0s}
\]
as long as $R-\gamma_0s>\de$. This proves the first part with $h_0:=h_\p$. The uniformity of $\gamma_0$ for small $r$ follows from the corresponding uniformity in
Proposition~\ref{evolution-phi-wulff-shape}.
\end{proof}

\subsection{Comparison with the evolution of $\psi$-Wulff shapes}\label{Comparison with the evolution of psi-Wulff shapes}
We now use Remark~\ref{Discr-avoid-remark} and Proposition~\ref{psi-Wulff-shape-evolution} to estimate the variation in time of the signed distance, as needed in the convergence
argument. As usual, we assume that $\psi$ is $\p$-regular.

Let $E^0$ be closed and let $(w,t)\in\R^N\times[0,T_h^*)$. Assume first that $d_h(w,t)>0$ and fix $\de\in(0,d_h(w,t))$, and let
\[
\de <R<d_h(w,t).
\]
Then $\dist^{\pso}(W_R^\psi(w),E_h(t))>0$. By Proposition~\ref{psi-Wulff-shape-evolution}, there exist $\gamma_0=\gamma_0(r,\de)>0$ and $h_0=h_0(r,\de)>0$ such that the discrete
evolution starting from $W_R^\psi(w)$ contains, after $m$ steps, $W_{R-\gamma_0mh}^\psi(w)$ as long as $mh<(R-\de)/\gamma_0$ and $0<h\leq h_0$. Between the times $t$ and $t+s$ at
most $m:=[s/h]+1$ discrete steps are performed. Hence, by the discrete avoidance principle, the evolved barrier $W_{R-\gamma_0mh}^\psi(w)$ stays in the complement of $E_h(t+s)$.
Therefore
\[
d_h(w,t+s)\ge R-\gamma_0mh\ge R-\gamma_0(s+h).
\]
Letting $R\nearrow d_h(w,t)$, we deduce that, for $0<h\le h_0$ and $0<s+h<(d_h (w, t)-\de)/\gamma_0$,
\beq\label{38}
    d_h(w,t+s)\ge d_h(w,t)-\gamma_0(s+h).
\eeq
If $d_h(w,t)<0$, we argue similarly.
For every $s$ with $0<s+h<(-d_h(w,t)-\de)/\gamma_0$ we may let $R\uparrow d_h(w,t)$ and obtain
\beq\label{39}
    d_h(w,t+s)\le d_h(w,t)+\gamma_0(s+h),
\eeq
for $0<h\le h_0$.

\subsection{Convergence of the scheme}
Throughout this section we again assume the $\p$-regularity of the mobility $\psi$. We prove the convergence, up to subsequences, of the minimising movement scheme introduced in
Section~\ref{min-mov-scheme-section} for an arbitrary closed initial set $E^0$. For each $h$ we fix any one of the discrete evolutions allowed by the construction above. By local
compactness for closed sets, there exist a sequence $h_l\downarrow0$, a closed set $E$, and an open set $A\subset E$ such that
\beq\label{28}
    \overline{E_{h_l}}\overset{\mathscr K}{\longrightarrow}E, \qquad \bigl(\mathrm{Int}(E_{h_l})\bigr)^c\overset{\mathscr K}{\longrightarrow}A^c,
\eeq
where closure, interior and complementation are taken relative to $\R^N\times[0,+\infty)$. For $t\ge0$ we set
\[
E(t):=\{w:(w,t)\in E\}, \qquad A(t):=\{w:(w,t)\in A\}.
\]

The following proposition establishes the compactness of the signed distance functions.

\begin{proposition}\label{convergence-Prop}
Possibly extracting a further subsequence (not relabelled), there exists a countable set $\mathcal N\subset(0,+\infty)$ such that, defining
\[
\td(w,t):= \dist^{\pso}(w,E(t)) - \dist^{\pso}(w,A^c(t)) \in\overline{\R},
\]
we have
\[
d_{h_l}(\cdot,t)\longrightarrow\td(\cdot,t) \qquad\text{locally uniformly in }\R^N
\]
for every $t\in(0,+\infty)\setminus\mathcal N$. 
Moreover, for every $w\in\R^N$, both functions $t\mapsto\dist^{\pso}(w,E(t))$ and $t\mapsto\dist^{\pso}(w,A^c(t))$ are left continuous and right lower semicontinuous. Finally,
$E(0)=E^0$ and $A(0)=\mathrm{Int}_{\R^N}(E^0)$.
\end{proposition}

\begin{proof}
The proof follows by the same compactness argument as in \cite[Proposition~4.4]{CMP17}, with the metric induced by $\pso$ and with \eqref{38} and \eqref{39} in place of the
corresponding time estimates.
\end{proof}

\begin{remark}\label{inequalities-distance-fun-remark}
Passing to the limit in \eqref{38} and \eqref{39}, we obtain the following one-sided estimates. If $(w,t)\in E^c$, $t\notin\mathcal N$, and $0<\de<\td(w,t)$, then there exists
$\gamma_0=\gamma_0(r,\de)>0$ such that
\beq\label{40}
    \td(w,t+s)\ge\td(w,t)-\gamma_0s
\eeq
for $0<s<(\td(w,t)-\de)/\gamma_0$. Analogously, if $(w,t)\in A$, $t\notin\mathcal N$, and $0<\de<-\td(w,t)$, then there exists $\gamma_0=\gamma_0(r,\de)>0$ such that
\beq\label{41}
    \td(w,t+s)\le\td(w,t)+\gamma_0s
\eeq
for $0<s<(-\td(w,t)-\de)/\gamma_0$. In particular, \eqref{40} implies property (b) of Definition~\ref{super-flow-def} for $E$, while \eqref{41} gives the analogous property for
$A^c$.
\end{remark}

\begin{remark}
As already observed in Remark~\ref{exitiction-time-remark}, if $E(t)=\emptyset$ for some $t\ge0$, then $E(s)=\emptyset$ for all $s\ge t$; hence we may define the extinction time
$T^{\mathrm{ext}}$ of $E$. Similarly, \eqref{41} implies that if $A(t)=\R^N$ for some $t\ge0$, then $A(s)=\R^N$ for all $s>t$, and we denote by $T^{\mathrm{all}}$ the extinction
time of $A^c$.
\end{remark}

Set $T^*:=\min\{T^{\mathrm{ext}},T^{\mathrm{all}}\}$. Then
\[
T^*\le\liminf_{l\to+\infty}T^*_{h_l}.
\]

The following theorem proves the convergence of the minimising movement scheme. The proof follows the argument of \cite[Theorem~3.20]{CMR26}. For the reader's convenience, we give
the details and explain where the estimates for the anisotropic Wulff shapes are used.

\begin{theorem}\label{convergence-min-mov-theo}
Let $E^0\subset\R^N$ be any closed set, and let $E$ and $A$ be the sets given in \eqref{28} for a family of discrete evolutions chosen as above. Assume that $\p$ is $\p$-regular.
Then $E$ is a weak superflow in the sense of Definition~\ref{super-flow-def} with initial datum $E^0$, while $A$ is a weak subflow with initial datum $E^0$.
\end{theorem}

\begin{proof}
Properties (a), (b) and (c) of Definition~\ref{super-flow-def} follow from Proposition~\ref{convergence-Prop} and Remark~\ref{inequalities-distance-fun-remark}; it remains to prove
(d). If $T^*=T^{\mathrm{all}}<T^{\mathrm{ext}}$, then $A(t)=\R^N$ and hence $E(t)=\R^N$ for $t\ge T^*$. Thus it is enough to verify (d) in $(0,T^*)$

Up to a further subsequence, and setting $z_{h_l}(\cdot,t):=0$ for $t\in(0,T^*)\setminus[h_l,T^*_{h_l})$, we may assume
\[
z_{h_l}\overset{*}{\wto}z \quad\text{in } L^\infty_{\mathbf w}\bigl(\R^N\times(0,T^*);\mathcal M^+(C_r^\p)\bigr), \qquad \|z\|_{L^\infty_{\mathbf w}}\le1.
\]
By \eqref{min mov scheme inequality1} and \eqref{min mov scheme},
\beq\label{61}
 h\divr_\p z_h^{k+1} \le d_{E_h^{k+1}}^{\pso}-d_{E_h^k}^{\pso} \qquad\text{in }\mathcal M(\{d_{E_h^{k+1}}^{\pso}>0\}),
\eeq
for all $l, k$ such that $k h_l < T^*_{h_l}$.
Let $\vphi\in C_c^\infty((\R^N\times(0,T^*))\setminus E)$ be non-negative. By the Kuratowski convergence, for $l$ large the support of $\vphi$ stays at positive distance from
$E_{h_l}$, and \eqref{61} yields
\[
\begin{split}
0\le \int_0^{T^*}\!\int_{\R^N}\Bigg[ &\vphi(w,t)\frac{d_{h_l}(w,t+h_l)-d_{h_l}(w,t)}{h_l}\\
&+\frac1{2r}\int_{C_r^\p} \bigl(\vphi(w+x,t)-\vphi(w+y,t)\bigr) \,d(z_{h_l})_{(w,t+h_l)}(x,y) \Bigg]dw\,dt\\
=\int_0^{T^*}\!\int_{\R^N}\Bigg[ &-\frac{\vphi(w,t)-\vphi(w,t-h_l)}{h_l}\,d_{h_l}(w,t)\\
&+\frac1{2r}\int_{C_r^\p} \bigl(\vphi(w+x,t)-\vphi(w+y,t)\bigr) \,d(z_{h_l})_{(w,t+h_l)}(x,y) \Bigg]dw\,dt .
\end{split}
\]
Passing to the limit by Proposition~\ref{convergence-Prop} and the weak-* convergence of $z_{h_l}$ gives
\[
\int_0^{T^*}\!\int_{\R^N}\left[ -\pa_t\vphi(w,t)\,\td(w,t) +\frac1{2r}\int_{C_r^\p} \bigl(\vphi(w+x,t)-\vphi(w+y,t)\bigr)\,dz_{(w,t)}(x,y) \right]dw\,dt\ge0.
\]
Therefore
\[
\pa_t\td\ge\divr_\p z \qquad\text{as measures in }(\R^N\times(0,T^*))\setminus E.
\]
This proves (d1) for $E$.

We next prove the local upper bound on $\divr_\p z$. Fix $\lambda>0$ and let $w$ satisfy $S:=d_{E_h^k}^{\pso}(w)\ge\lambda$. Choose $\xi\in E_h^k$ with $\pso(w-\xi)=S$ and set
$S_0:=\e_0S$. By $\p$-regularity there exists $\widetilde w$ such that
\[
W_{S_0}^\p(\widetilde w)\subset W_S^\psi(\xi), \qquad w\in\pa W_{S_0}^\p(\widetilde w).
\]
Since signed distances are monotone with respect to inclusion, \eqref{44} gives, for every $y$,
\beq\label{anisotropic distance upper barrier}
\begin{split}
d_{E_h^k}^{\pso}(y) &\le \pso(y-\xi) =S+d_{W_S^\psi(\xi)}^{\pso}(y)\\
&\le S+d_{W_{S_0}^\p(\widetilde w)}^{\pso}(y)\\
&\le S+c_2(\po(y-\widetilde w)-S_0)_+ -c_1(\po(y-\widetilde w)-S_0)_-.
\end{split}
\eeq
Set $\eta_\lambda:=-c_1\e_0\lambda/2$. Since $S_0\ge\e_0\lambda$, Proposition~\ref{Wulff-shape-comparison}, together with the monotonicity of its constants with respect to the
radius, gives constants $c_\lambda=c_\lambda(r)>0$ and $h_\lambda=h_\lambda(r)>0$, independent of $S$. Recall also that $c_\lambda$ can be chosen independent of $r$ for $r$
sufficiently small. The right-hand side of \eqref{anisotropic distance upper barrier} is a datum in $\mathscr L^+$, so the one-sided comparison in Proposition~\ref{comparison}
applies even when $d_{E_h^k}^{\pso}$ is noncoercive and yields for any $y$
\[
u_h^{k+1}(y)\leq S+\big[c_2(\po(y-\widetilde w)-S_0)_+ -c_1(\po(y-\widetilde w)-S_0)_-+c_\lambda h \big]\lor\eta_\lambda.
\]
Evaluating at $y=w$ and recalling that $\po(w-\widetilde w)=S_0$ gives
\beq\label{56}
    u_h^{k+1}(w)\le d_{E_h^k}^{\pso}(w)+c_\lambda h
\eeq
whenever $d_{E_h^k}^{\pso}(w)\ge\lambda$ and $h\le h_\lambda$. By \eqref{min mov scheme},
\[
\divr_\p z_h^{k+1}\le c_\lambda \qquad\text{as measures on }\{d_{E_h^k}^{\pso}\ge\lambda\}.
\]
Passing to the limit, we obtain
\[
\divr_\p z\le c_{\lambda} \qquad\text{as measures on } \{(w,t):0<t<T^*,\ \td(w,t)>\lambda\}.
\]
This proves (d3).

It remains to identify the limiting Cahn-Hoffman field. If $S=d_{E_h^k}^{\pso}(w)>0$, the $\pso$-Lipschitz property of the signed distance and \eqref{costants psi phi} give
\[
d_{E_h^k}^{\pso}(y) \ge S-\pso(y-w) \ge S-c_2\po(y-w).
\]
Recall that $\theta^\tau$ denotes the solution with datum $\po$ and time step $\tau$, see Proposition~\ref{comparison-phi-Wulff-shape}. If $z^\tau$ is a corresponding Cahn-Hoffman field, then
\[
-\tau\divr_\p z^\tau+\theta^\tau=\po.
\]
After translation, multiplication by $-c_2$, and exchange of the two coordinates in $z^\tau$, the function $S-c_2\theta^\tau(\cdot-w)$ is therefore a solution, with time step
$c_2\tau$, corresponding to the datum $S-c_2\po(\cdot-w)$. Choosing $\tau=h/c_2$, and using
\[
d_{E_h^k}^{\pso} (\cdot) \ge S-c_2\po(\cdot-w),
\]
the comparison principle yields
\[
u_h^{k+1}(w)\ge S-c_2\theta^{h/c_2}(0).
\]
Taking $\alpha=1/(2(N+1))$ in \eqref{estimate zero Wulff shape}, and combining this inequality with \eqref{56}, we obtain, for every $\lambda>0$,
\beq\label{discrete u distance convergence}
\|u_{h_l}(\cdot,t)-d_{h_l}(\cdot,t-h_l)\|_{L^\infty(\{d_{h_l}(\cdot,t-h_l)\ge\lambda\})} \le C_\lambda h_l^{1/(2(N+1))}
\eeq
for $l$ large and $t\in[h_l,T^*_{h_l})$. The analogous estimate on $\{d_{h_l}(\cdot,t-h_l)\le-\lambda\}$ follows by applying the same one-sided comparisons to $-u_{h_l}$ and
$-d_{h_l}$, with the two coordinates in the Cahn-Hoffman fields exchanged. By the Helly-type compactness argument used in the proof of Proposition~\ref{convergence-Prop} (see also
\cite[formula~(4.13)]{CMP17}), for a.e. $t\in(0,T^*)$ we have $d_{h_l}(\cdot,t-h_l)\to\td(\cdot,t)$ locally uniformly. Moreover,
\eqref{min mov scheme inequality1}-\eqref{min mov scheme inequality2} give $|u_{h_l}(\cdot,t)|\le |d_{h_l}(\cdot,t)|$. Hence, on $\{|d_{h_l}(\cdot,t-h_l)|<\lambda\}$,
\[
|u_{h_l}(\cdot,t)-d_{h_l}(\cdot,t-h_l)| \le 2\lambda+o(1)
\]
locally uniformly. Combining this with the preceding estimates gives
\beq\label{57}
    u_{h_l}(\cdot,t)\longrightarrow\td(\cdot,t) \qquad\text{locally uniformly in }\R^N.
\eeq
 The same compactness argument and the inequality $|u_{h_l}|\le|d_{h_l}|$ also yield local uniform boundedness of $(u_{h_l})_l$ in space-time.

We finally pass to the limit in the Cahn-Hoffman identity. Fix $\zeta\in C_c^\infty(\R^N\times(0,T^*))$. For $l$ sufficiently large, the time projection of
$\operatorname{supp}\zeta$ is contained in $[h_l,T^*_{h_l})$, and the discrete Cahn-Hoffman identity gives
\[
\begin{split}
&\int_0^{T^*}\!\int_{\R^N}\zeta(w,t) \int_{C_r^\p}\bigl(u_{h_l}(w+x,t)-u_{h_l}(w+y,t)\bigr) \,d(z_{h_l})_{(w,t)}(x,y)\,dw\,dt\\
&\qquad=\int_0^{T^*}\!\int_{\R^N} \zeta(w,t)\osc_{W_r^\p(w)}u_{h_l}(\cdot,t)\,dw\,dt.
\end{split}
\]
By \eqref{57}, local boundedness and dominated convergence, the right-hand side converges to the same expression with $\td$ in place of $u_{h_l}$. For the left-hand side, add and
subtract the same term with $u_{h_l}$ replaced by $\td$, using \eqref{57}, the mass bound, and the weak-* convergence of $z_{h_l}$ to conclude that
\begin{multline*}
\int_0^{T^*}\!\int_{\R^N}\zeta(w,t) \int_{C_r^\p}\bigl(u_{h_l}(w+x,t)-u_{h_l}(w+y,t)\bigr) \,d(z_{h_l})_{(w,t)}(x,y)\,dw\,dt \\
\to \int_0^{T^*}\!\int_{\R^N}\zeta(w,t) \int_{C_r^\p}\bigl(\td(w+x,t)-\td(w+y,t)\bigr) \,dz_{(w,t)}(x,y)\,dw\,dt
\end{multline*}
Hence, we obtain
\[
\int_{C_r^\p} \bigl(\td(w+x,t)-\td(w+y,t)\bigr)\,dz_{(w,t)}(x,y) = \osc_{W_r^\p(w)}\td(\cdot,t)
\]
for a.e. $(w,t)$. By Remark~\ref{Remark-Cahn-Hoffman}, the same field is also a generalised Cahn-Hoffman field for the positive part $\td^+=\dist^{\pso}(\cdot,E(t))$. Indeed,
taking the positive part is a non-decreasing operation, so the pairs attaining the maximum and minimum in \eqref{27} still attain the maximum and minimum after truncation. This
proves (d2), and therefore $E$ is a weak superflow.

The proof that $A$ is a weak subflow is completely analogous. More precisely, one applies the preceding argument to the negative part of the signed distance, replacing
$(u_h,d_h,z_h)$ with $(-u_h,-d_h,\check z_h)$, where $\check z_h$ is obtained by exchanging the two coordinates in $z_h$. Then $\divr_\p\check z_h=-\divr_\p z_h$, and the previous
estimates are unchanged after this transformation. Hence $A^c$ is a weak superflow with initial datum $\R^N\setminus\operatorname{Int}_{\R^N}(E^0)$, and therefore $A$ is a weak
subflow with initial datum $E^0$.
\end{proof}
Motivated by the above theorem, we conclude the section with the following definition.
\begin{definition}\label{def:variationalflows}
Let $E^0\subset\R^N$ be closed. A weak superflow $E$ of \eqref{curvature} with initial datum $E^0$ is called a \textit{variational superflow} if it is obtained as a Kuratowski
limit of minimising movements, namely if there exist $h_l\downarrow0$ and discrete evolutions $E_{h_l}$ constructed according to the minimising movement scheme starting from $E^0$
such that
\[
\overline{E_{h_l}}\overset{\mathscr K}{\longrightarrow}E.
\]
Similarly, a weak subflow $A$ of \eqref{curvature} with initial datum $E^0$ is called a \textit{variational subflow} if there exist $h_l\downarrow0$ and corresponding minimising
movements $E_{h_l}$ such that
\[
\bigl(\operatorname{Int}E_{h_l}\bigr)^c \overset{\mathscr K}{\longrightarrow}A^c.
\]
Finally, $E$ is a {\em variational weak flow} of \eqref{curvature} starting from $E_0$ if it is a variational weak superflow and $\operatorname{Int}E$ is a variational weak subflow
of \eqref{curvature} starting from $E_0$. Here, the interior is taken in the topology of $\R^N\times[0,+\infty)$.
\end{definition}

\subsection{Level set flows with compact boundaries} We start by stating the following comparison principle.
\begin{theorem}\label{comparison-principle-sets}
    Let $E^0, F^0 \subset \R^N$ be two closed sets, and assume that at least one of $E^0$ and $F^0$ has compact boundary. Let $E,F$ be, respectively, a superflow with initial datum
    $E^0$ and a subflow with initial datum $F^0$, according to Definition \ref{super-flow-def}. If $\dist^{\pso}(E^0, \R^N \setminus F^0) \eqqcolon \Delta >0$, then
    \[
    \dist^{\pso}(E(t), \R^N \setminus F(t)) \geq \Delta \hspace{5mm} \text{for every } t \geq 0,
    \]
    with the convention that $\dist^{\pso}(A, \emptyset) = + \infty$ for every set $A$.
\end{theorem}
\begin{proof}
 The proof is the same as that of \cite[Theorem~3.21]{CMR26}. Note that the argument only requires one of the two sets to have compact boundary.
\end{proof}

An immediate corollary of Theorem \ref{comparison-principle-sets} is the following comparison principle between level set sub- and super-solutions.

\begin{theorem}
    Let $u_0,v_0$ be uniformly continuous functions on $\R^N$ such that $\{ u_0 \leq \lambda \}, \{ v_0 \leq \lambda \}$ have compact boundary for every $\lambda \in \R$. Let $u,v$
    be respectively a level set sub-solution with initial datum $u_0$ and a level set super-solution with initial datum $v_0$, according to Definition \ref{sup-sub-levle-set-def}.
    If $u_0 \leq v_0$, then $u \leq v$.
\end{theorem}

\begin{proof}
    The proof follows by arguing exactly as in \cite[Theorem 2.8]{CMP17}.
\end{proof}

We are now in a position to prove existence and uniqueness of the level set flow by standard arguments. We introduce the following $\textit{level set minimizing movements}$.

Let $u_0 \colon \R^N \rightarrow \R$ be a uniformly continuous function such that $\{u_0 \leq \lambda \}$ has compact boundary for every $\lambda \in \R$. Let
$E^{0,\lambda} \coloneqq \{u_0 \leq \lambda \}$ and let $E_{\lambda,h}$ be the corresponding time discrete evolutions, defined according to \eqref{min mov scheme} with $E^0$
replaced by $E^{0,\lambda}$.

We introduce the level set discrete evolution $w_h \colon \R^N \times \R \rightarrow \R$ defined by
\beq\label{70}
    w_h(x,t) \coloneqq \inf \{ \lambda \in \R \, : \, x \in E_{\lambda, h}(t)\}.
\eeq

Using Proposition \ref{evolution-phi-wulff-shape} and arguing as in \cite[Section~4.3]{CMNP19} and \cite[Lemma 6.13]{CMP12}, we obtain the following equicontinuity lemma.
\begin{lemma}\label{equi-continuity-lemma}
    Let $w_h$ be as in \eqref{70}. Then, for every $h>0$ and for every $t>0$, the function $w_h( \cdot, t)$ has the same modulus of continuity as $u_0$. Moreover, for every
    $\e >0$, there exist $\tau = \tau (r,\e) >0$ and $h_0 = h_0 (r, \e)>0$ such that for every $|t-t'| \leq \tau$ and $h \leq h_0$ we have $|w_h (\cdot,t)- w_h (\cdot, t')| < \e$.
    In addition, for sufficiently small values of $r$, the constant $\tau$ can be chosen independent of $r$, that is $\tau = \tau(\e)$.
\end{lemma}

Arguing as in the proof of \cite[Theorem 4.8]{CMNP19}, we finally obtain the following main result on the well-posedness of the level set flow.

\begin{theorem}\label{level-set-flow-theorem}
    Assume that $\psi$ is $\p$-regular. Let $u_0$ be a uniformly continuous function on $\R^N$ such that $\{ u_0 \leq \lambda \}$ has compact boundary for every $\lambda \in \R$.
    Then the following holds:
    \begin{itemize}
        \item[(i)] $(\mathrm{Existence \, and \, uniqueness})$ There exists a unique solution $u$ to the level set flow with initial datum $u_0$, in the sense of Definition
        \ref{sup-sub-levle-set-def}.
        \item[(ii)] $\mathrm{(Approximation \, via \, minimizing \, movements)}$. The solution $u$ is the locally uniform limit in $\R^N \times [0,+ \infty)$ as $h \rightarrow 0^+$
        of the level set minimising movements $w_h$ defined in \eqref{70}.
        \item[(iii)] $\mathrm{(Properties \, of \, the \, level \, set \, flow)}$. For all but countably many $\lambda \in \R$, the fattening phenomenon does not occur; that is,
        \beq
        \begin{array}{rcl}\label{71}
\{(x,t)\,:\, u(x, t) < \lambda\}& = &\mathrm{Int\,}(\{(x,t)\,:\, u(x, t) \le \lambda\})\,, \vspace{5pt} \\
\overline {\{(x,t)\,:\, u(x, t) < \lambda\}} &=& \{(x,t)\,:\, u(x, t) \le \lambda\}\,.
\end{array}
\eeq
        Moreover, for every $\lambda$ such that \eqref{71} holds, the sublevel set $\{ u \leq \lambda\}$ is the unique solution to \eqref{curvature} in the sense of Definition
        \ref{super-flow-def}, with initial datum $E^{0,\lambda}$, and
       \beq
            E_{\lambda,h} \overset{\mathscr{K}}{\rightarrow } \{ (x,t) \, : \, u(x,t) \leq \lambda \} \hspace{3mm} \text{and} \hspace{3mm} (\mathrm{Int} E_{\lambda, h})^c
            \overset{\mathscr{K}}{\rightarrow} \{(x,t) \, : \, u(x,t) \geq \lambda \}
        \eeq
        as $h \rightarrow 0^+$. Finally, for every $\lambda \in \R$ the sets $\{ (x,t) \, : \, u(x,t) \leq \lambda \}$ and $\{ (x,t) \, : \, u(x,t) < \lambda \}$ are, respectively,
        the maximal superflow and minimal subflow with initial datum $E^{0,\lambda}$.
    \end{itemize}
\end{theorem}
In the following, whenever $\lambda$ is such that \eqref{71} holds, we say that the set $E^{0, \lambda}$ is non-fattening for the flow. We will use the same terminology also for the local anisotropic flow \eqref{classical curvature}.

\section{Asymptotic behaviour as $r \rightarrow 0^+$}\label{sec:asymptotic}

In this section we study the local limit of the anisotropic Minkowski flows as the interaction radius tends to zero. As usual, we assume that $\psi$ is $\p$-regular. The proof
consists of two main steps. First, the estimates for the evolution of Wulff shapes proved in the previous section yield compactness of the distance functions, uniformly for small
$r$. We then average the measure-valued nonlocal Cahn-Hoffman fields at scale $r$ and show that the resulting vector fields converge to a field taking values in $\pa\p(\nabla d)$.
For the latter step, the only additional assumption is a uniform upper bound on the nonlocal divergences away from the evolving interface.

Let $E^0\subset\R^N$ be closed and, for any $r>0$, let $E_r$ be a weak superflow of \eqref{curvature} with initial datum $E^0$, in the sense of Definition~\ref{super-flow-def}. We
set
\[
T_r^{\mathrm{ext}}:=\inf\{t>0:E_r(s)=\emptyset\text{ for every }s\ge t\}, \qquad d_r(w,t):=\dist^{\pso}(w,E_r(t))
\]
for $0<t<T_r^{\mathrm{ext}}$, and we choose an associated generalised Cahn-Hoffman field for $d_r(\cdot,t)$ at scale $r$, for a.e. $t$,
\[
z^r\in L^\infty_{\mathbf w}\bigl(\R^N\times(0,T_r^{\mathrm{ext}});\mathcal M^+(C_r^\p)\bigr)
\]
as in Definition~\ref{super-flow-def}.

We say that the family $(E_r)_{r>0}$ satisfies condition $(\ND)$ if the positive part of the corresponding nonlocal divergences is bounded uniformly in $r$ on sets which stay a
positive distance from the interface: for every $\de>0$ there is $c_{\de}>0$, independent of $r$ for $r$ small, such that
\beq\label{bound divergence}
\divr_\p z^r\le c_{\de} \quad\text{as measures on } \{(w,t):0<t< T_r^{\mathrm{ext}},\ d_r(w,t)\ge\de\}.
\eeq

The main result of the section is the following.

\begin{theorem}\label{asymptotic-behavior-theorem}
Assume $\psi$ to be $\p$-regular. Let $E^0\subset\R^N$ be closed and, for every $r>0$, let $E_r$ be a variational weak superflow of \eqref{curvature} with initial datum $E^0$, see
Definition~\ref{def:variationalflows}. Then there exist $r_k\downarrow0$ and a closed set $E\subset\R^N\times[0,+\infty)$ such that
\[
E_{r_k}\overset{\mathscr K}{\longrightarrow}E,
\]
and $E$ is a weak superflow, with initial datum $E^0$, of the local anisotropic flow
\beq\label{classical curvature}
V=-\psi(\nu)H^\p
\eeq
in the sense of Definition~\ref{classic-superflow-def}. More generally, the same conclusion holds if $(E_r)$ is a family of weak superflows starting from $E^0$ and satisfying
condition $(\ND)$.
\end{theorem}

We first state two consequences that will be used later.

\begin{corollary}\label{weak-asymptotic-conv-cor-1}
Assume $\psi$ to be $\p$-regular. Let $E^0=\overline{\operatorname{Int}(E^0)}$ be non-fattening for \eqref{classical curvature}. For every $r>0$, let $A_r\subset E_r$ be a
variational weak subflow and superflow of \eqref{curvature}, respectively, with initial datum $E^0$. Then
\[
E_r\overset{\mathscr K}{\longrightarrow}E \qquad\text{as }r\downarrow0,
\]
where $E$ is the unique weak flow of \eqref{classical curvature} with initial datum $E^0$. In particular, the same conclusion holds if for every $r>0$, $E_r$ is a variational weak
flow of \eqref{curvature} starting from $E^0$.
\end{corollary}

\begin{proof}
Apply Theorem~\ref{asymptotic-behavior-theorem} simultaneously to $E_r$ and to $(A_r)^c$. Along every sequence $r_k\downarrow0$ one can therefore extract a subsequence, not
relabelled, such that
\[
E_{r_k}\overset{\mathscr K}{\longrightarrow}E, \qquad A_{r_k}^c\overset{\mathscr K}{\longrightarrow}A^c,
\]
with $E$ a weak superflow and $A$ a local weak subflow of \eqref{classical curvature}. Let $u$ be the local level set flow of \eqref{classical curvature} with initial datum given
by the signed $\pso$-distance from $E^0$. By the comparison theory of \cite{CMP17,CMNP19},
\[
\{u<0\}\subset A\subset E\subset\{u\le0\}.
\]
The non-fattening assumption (see \eqref{71}) gives $\overline{\{u<0\}}=\{u\le0\}$, hence $E=\overline A=\{u\le0\}$. Thus every sequence of radii tending to zero has the same
Kuratowski limit, and the whole family converges.
\end{proof}

\begin{corollary}\label{level-set-asymptotic-cor}
Let $\psi$ be $\phi$-regular. Let $u^0:\R^N\to\R$ be uniformly continuous and assume that $\{u^0\le\lambda\}$ has compact boundary for every $\lambda\in\R$. For $r>0$, let $u_r$ be
the level set flow of \eqref{curvature} constructed in Theorem~\ref{level-set-flow-theorem}. Then
\[
u_r\longrightarrow u \quad\text{locally uniformly in }\R^N\times[0,+\infty)
\]
as $r\downarrow0$, where $u$ is the unique level set flow of \eqref{classical curvature} with initial datum $u^0$.
\end{corollary}

\begin{proof}
The functions $u_r$ have the same spatial modulus of continuity as $u^0$. The time modulus obtained from the Wulff barriers in Lemma~\ref{equi-continuity-lemma} can be chosen
uniformly for all sufficiently small $r$, by the last assertion of Proposition~\ref{psi-Wulff-shape-evolution}. Hence $(u_r)_{r>0}$ is locally equicontinuous and every sequence
$r_k\downarrow0$ admits a subsequence, not relabelled, and a function $u:\R^N\times[0,+\infty)\to\R$ such that
\[
u_{r_k}\longrightarrow u \qquad\text{locally uniformly in }\R^N\times[0,+\infty).
\]
In particular, $u(\cdot,0)=u^0$.

There is a common set $\Lambda\subset\R$ of full measure such that, for every level $\lambda$ in this set, the closed sets $E_k^\lambda:=\{u_{r_k}\le\lambda\}$ and
$F_k^\lambda:=\{u_{r_k}\ge\lambda\}$ are variational weak superflows of \eqref{curvature}, with initial data $\{u^0\le\lambda\}$ and $\{u^0\ge\lambda\}$, respectively. This follows
directly from the minimising movement construction in Theorem~\ref{level-set-flow-theorem}.

Shrinking $\Lambda$, if necessary, we may also assume that $\overline{\{u<\lambda\}}=\{u\le\lambda\}$ and $\overline{\{u>\lambda\}}=\{u\ge\lambda\}$ for every $\lambda\in\Lambda$
(non-fattening condition). For all such $\lambda$, the local uniform convergence implies
\[
E_k^\lambda\overset{\mathscr K}{\longrightarrow}\{u\le\lambda\}, \qquad F_k^\lambda\overset{\mathscr K}{\longrightarrow}\{u\ge\lambda\}.
\]
Indeed, if $x_k\in E_k^\lambda$ and $x_k\to x$, then $u(x)\le\lambda$, while every $x\in\{u<\lambda\}$ belongs to $E_k^\lambda$ for all sufficiently large $k$. The first
convergence follows from $\overline{\{u<\lambda\}}=\{u\le\lambda\}$, and the second one is proved analogously.

We can therefore apply Theorem~\ref{asymptotic-behavior-theorem} to the two variational families $(E_k^\lambda)_k$ and $(F_k^\lambda)_k$. It follows that, for every
$\lambda\in\Lambda$, $\{u\le\lambda\}$ is a weak superflow of \eqref{classical curvature} with initial datum $\{u^0\le\lambda\}$, while $\{u\ge\lambda\}$ is a weak superflow of
\eqref{classical curvature} with initial datum $\{u^0\ge\lambda\}$. Thus $u$ is both a level set supersolution and a level set subsolution of \eqref{classical curvature}, and hence
a level set solution.

By uniqueness of the local level set solution, $u$ is independent of the chosen subsequence. Therefore $u_r\to u$ locally uniformly in $\R^N\times[0,+\infty)$ as $r\downarrow0$.
\end{proof}

\subsection*{Compactness of the distance functions}

Fix a sequence $r_k\downarrow0$. By the local compactness of closed sets for Kuratowski convergence, after extracting a subsequence we may assume
\beq\label{kuratowski convergence asymptotic}
E_{r_k}\overset{\mathscr K}{\longrightarrow}E
\eeq
for some closed $E\subset\R^N\times[0,+\infty)$.

The Wulff barriers give the uniform one-sided time estimate needed below. Let $0<\de<d_{r_k}(w,t)$ and choose $R$ with
\[
\de<R<d_{r_k}(w,t).
\]
Then $W_R^\psi(w)$ and $E_{r_k}(t)$ are separated by
\[
\Delta:=\dist^{\pso}\bigl(E_{r_k}(t),W_R^\psi(w)\bigr) =d_{r_k}(w,t)-R>0.
\]
Let $F$ be the weak superflow starting from $W_R^\psi(w)$, obtained as a limit of the minimising movement scheme. By the comparison principle of
Theorem~\ref{comparison-principle-sets}, the separation $\Delta$ is preserved, while Proposition~\ref{psi-Wulff-shape-evolution} gives, for a constant $\gamma_\de>0$ independent of
all sufficiently small $r_k$,
\[
W_{R-\gamma_\de s}^\psi(w)\subset F(s) \qquad\text{for }0<s<(R-\de)/\gamma_\de.
\]
Consequently
\[
d_{r_k}(w,t+s) \ge \Delta+R-\gamma_\de s =d_{r_k}(w,t)-\gamma_\de s.
\]
Thus, whenever $d_{r_k}(w,t)>\de$,
\beq\label{72}
d_{r_k}(w,t+s)\ge d_{r_k}(w,t)-\gamma_\de s
\eeq
for $s>0$ sufficiently small, with $\gamma_\de$ independent of $k$.

The Helly-type compactness argument of \cite[Proposition~4.4]{CMP17}, with \eqref{72} in place of the corresponding discrete estimate, now yields the following.

\begin{proposition}\label{distance-function-convergence}
Possibly extracting a further subsequence, there is a countable set $\mathcal N\subset(0,+\infty)$ such that, for every $t\notin\mathcal N$,
\beq\label{distance local uniform}
d_{r_k}(\cdot,t)\longrightarrow d(\cdot,t):=\dist^{\pso}(\cdot,E(t)) \quad\text{locally uniformly in }\R^N,
\eeq
with locally uniform bounds in space-time. Moreover, for every $w\in\R^N$, $t\mapsto\dist^{\pso}(w,E(t))$ is left continuous and right lower semicontinuous, and $E(0)\subset E^0$.
\end{proposition}


\subsection*{Limit of the nonlocal Cahn-Hoffman fields}

Extend $z^{r_k}$ by zero for $t\ge T_{r_k}^{\mathrm{ext}}$ and define
\beq\label{z hat definitions}
\widehat z_{r_k}(w,t) :=\frac1{2r_k}\int_{C_{r_k}^\p}(x-y)\,dz^{r_k}_{(w,t)}(x,y).
\eeq
Since $\po(x-y)\le2r_k$ for $(x, y) \in C_{r_k}^\p$ and $z^{r_k}_{(w,t)}(C_{r_k}^\p)\le1$,
\beq\label{uniform bound norm z hat}
\po(\widehat z_{r_k}) \le\frac1{2r_k}\int_{C_{r_k}^\p}\po(x-y)\,dz^{r_k}_{(w,t)}(x,y) \le1
\eeq
for a.e. $(w,t)$. Hence, after extracting a further subsequence,
\beq\label{weak convergence hat z}
\widehat z_{r_k}\overset{*}{\wto}\widehat z \quad\text{in }L^\infty_{\mathrm{loc}}(\R^N\times(0,+\infty);\R^N),
\eeq
and $\po(\widehat z)\le1$ a.e.

In the following proposition we show that the nonlocal divergences converge to the ordinary divergence of $\widehat z$ and prove the upper bound required in the definition of a
local weak superflow.

\begin{proposition}\label{divergence-mean-inequality-measures-prop}
Let $d=\dist^{\pso}(\cdot,E(\cdot))$. Then
\beq\label{divergence mean inequality measures}
\pa_t d\ge\Div\widehat z \quad\text{in }(\R^N\times(0,T^{\mathrm{ext}}))\setminus E
\eeq
in the sense of distributions, where $T^{\mathrm{ext}}$ is the extinction time of $E$. Moreover, if condition $(\ND)$ holds, then for every $\de>0$ there exists $C_{\de}>0$ such
that
\beq\label{local div limit bound}
\Div\widehat z\le C_{\de} \quad\text{in }\{(w,t):0<t<T^{\mathrm{ext}},\ d(w,t)\ge\de\}.
\eeq
\end{proposition}

\begin{proof}
Let $\eta\in C_c^\infty((\R^N\times(0,T^{\mathrm{ext}}))\setminus E)$. For $k$ large, the support of $\eta$ lies before $T_{r_k}^{\mathrm{ext}}$ and is separated from $E_{r_k}$. 
By definition of the nonlocal divergence,
\[
\langle\mathfrak{Div}_\p^{r_k}z^{r_k},\eta\rangle =-\frac1{2r_k}\int\!\int_{C_{r_k}^\p} \bigl(\eta(w+x,t)-\eta(w+y,t)\bigr)\,dz^{r_k}_{(w,t)}(x,y)\,dw\,dt.
\]
Taylor's formula gives, uniformly on the relevant compact set,
\[
\eta(w+x,t)-\eta(w+y,t) =\nabla\eta(w,t)\cdot(x-y)+R_k(w,t,x,y),
\]
with $|R_k|\le C r_k^2$. Since the mass of $z^{r_k}_{(w,t)}$ is at most one, we have
\beq\label{div convergence test}
\langle\mathfrak{Div}_\p^{r_k}z^{r_k},\eta\rangle =-\int\!\int\nabla\eta\cdot\widehat z_{r_k}\,dw\,dt+O(r_k).
\eeq
Using \eqref{weak convergence hat z}, we obtain
\[
\mathfrak{Div}_\p^{r_k}z^{r_k}\longrightarrow\Div \, \widehat z \qquad\text{in }\mathcal D'((\R^N\times(0,T^{\mathrm{ext}}))\setminus E).
\]
Passing to the limit in the inequalities $\pa_t d_{r_k}\ge\mathfrak{Div}_\p^{r_k}z^{r_k}$, using Proposition~\ref{distance-function-convergence} and the uniform one-sided time
estimate \eqref{72}, gives \eqref{divergence mean inequality measures}.

Finally, let $\eta\ge0$ be supported in $\{0<t<T^{\mathrm{ext}},\ d\ge\de\}$. Then, for $k$ large its support is contained in $\{d_{r_k}\ge\de/2\}$; hence condition $(\ND)$ and
\eqref{div convergence test} give
\[
\langle\Div\widehat z,\eta\rangle =\lim_{k\to\infty}\langle\mathfrak{Div}_\p^{r_k}z^{r_k},\eta\rangle \le c_{\de/2}\int\eta.
\]
This proves \eqref{local div limit bound}.
\end{proof}

We next prove a lower semicontinuity result relating the nonlocal oscillations to the anisotropic total variation of the limiting distance function.

\begin{lemma}\label{osc-lemma}
For every non-negative $\eta\in C_c^\infty((\R^N\times(0,+\infty))\setminus E)$,
\beq\label{osc liminf}
\liminf_{k\to\infty}\frac1{2r_k}\int\eta(w,t) \osc_{W_{r_k}^\p(w)}d_{r_k}(\cdot,t)\,dw\,dt \ge\int\eta(w,t)\p(\nabla d(w,t))\,dw\,dt.
\eeq
\end{lemma}

\begin{proof}
LEt $\mathcal N$ be the set given by Proposition \ref{distance-function-convergence}. Fix $t\notin\mathcal N$. By \eqref{distance local uniform}, $d_{r_k}(\cdot,t)\to d(\cdot,t)$ in $L^1_{\mathrm{loc}}$. The $\Gamma$-liminf inequality for the anisotropic outer
Minkowski content, applied to $d_{r_k}$ and $-d_{r_k}$, gives (see \cite[Corollary~3.2]{CLL14})
\beq\label{osc open liminf}
\liminf_{k\to\infty}\frac1{2r_k}\int_U \osc_{W_{r_k}^\p(w)}d_{r_k}(\cdot,t)\,dw \ge\int_U\p(\nabla d(w,t))\,dw
\eeq
for every open $U\Subset\R^N$.

Set
\[
\mu_k^t:=\frac1{2r_k}\osc_{W_{r_k}^\p(\cdot)}d_{r_k}(\cdot,t)\,\mathcal L^N, \qquad \mu^t:=\p(\nabla d(\cdot,t))\,\mathcal L^N.
\]
Equation \eqref{osc open liminf} and integration over the superlevel sets imply, for every non-negative $\zeta\in C_c(\R^N)$,
\[
\liminf_{k\to\infty}\int\zeta\,d\mu_k^t\ge\int\zeta\,d\mu^t.
\]
Indeed, apply \eqref{osc open liminf} to the open sets $\{\zeta>s\}$ and then use Fatou's lemma in the variable $s$. Taking $\zeta=\eta(\cdot,t)$ and applying Fatou's lemma once
more in time proves \eqref{osc liminf}.
\end{proof}

It remains to show that the weak-* limit $\widehat z$ is a local Cahn-Hoffman field for $d$ away from $E$.

\begin{proposition}\label{Cahn-hoffman-proposition-asymptotic}
Under the assumptions of Theorem~\ref{asymptotic-behavior-theorem},
\beq\label{sub gradient zeta}
\widehat z\in\pa\p(\nabla d) \quad\text{a.e. in }(\R^N\times(0,T^{\mathrm{ext}}))\setminus E.
\eeq
\end{proposition}
\begin{proof}
By \eqref{uniform bound norm z hat}-\eqref{weak convergence hat z}, $\po(\widehat z)\le1$. In view of the characterization \eqref{sub gradient}, it remains to prove
\beq\label{CH limit target}
\widehat z\cdot\nabla d\ge\p(\nabla d) \qquad\text{a.e. outside }E;
\eeq
the reverse inequality follows immediately from the polar inequality.

Let $0\le\eta\in C_c^\infty((\R^N\times(0,T^{\mathrm{ext}})) \setminus E)$. By the Cahn-Hoffman identity \eqref{CH} and Lemma~\ref{osc-lemma},
\beq\label{CH liminf start}
\liminf_{k\to\infty}I_k \ge \int\eta\,\p(\nabla d),
\eeq
where
\[
I_k:=\frac1{2r_k}\int\eta(w,t) \int_{C_{r_k}^\p} \bigl(d_{r_k}(w+x,t)-d_{r_k}(w+y,t)\bigr) \,dz^{r_k}_{(w,t)}(x,y)\,dw\,dt.
\]
We claim that
\beq\label{CH limsup target}
\limsup_{k\to\infty}I_k \le\int\eta\,\widehat z\cdot\nabla d.
\eeq

Fix $\e>0$ and let $d_\e=\rho_\e*d$ be the spatial mollification of $d$. Choose compact sets $K\Subset\R^N$ and $J\Subset(0,T^{\mathrm{ext}})$ such that, for all sufficiently large
$k$, the support of $\eta$ and all the spatial translates occurring below are contained in $K\times J$. Set
\[
m_{\e,k}(t) := \max_{w\in K} \bigl(d_{r_k}(w,t)-d_\e(w,t)\bigr), \qquad q_{\e,k} := d_{r_k}-d_\e-m_{\e,k}.
\]
Thus $q_{\e,k}\le0$ on $K\times J$. Since $m_{\e,k}$ depends only on $t$, it cancels when taking spatial differences. Therefore
\[
I_k= \frac1{2r_k}\int\eta(w,t) \int_{C_{r_k}^\p} \bigl(d_\e(w+x,t)-d_\e(w+y,t)\bigr) \,dz^{r_k}_{(w,t)}(x,y)\,dw\,dt +R_{\e,k},
\]
where
\[
R_{\e,k} := \frac1{2r_k}\int\eta(w,t) \int_{C_{r_k}^\p} \bigl(q_{\e,k}(w+x,t) -q_{\e,k}(w+y,t)\bigr) \,dz^{r_k}_{(w,t)}(x,y)\,dw\,dt.
\]
Since $d_\e$ is smooth in space, Taylor's formula and $z^{r_k}_{(w,t)}(C_{r_k}^\p)\le1$ give, for fixed $\e$,
\[
\frac1{2r_k}\int\eta \int_{C_{r_k}^\p} \bigl(d_\e(w+x)-d_\e(w+y)\bigr)\,dz^{r_k}_{(w,t)} = \int\eta\,\nabla d_\e\cdot\widehat z_{r_k}+O_\e(r_k).
\]
Hence, by \eqref{weak convergence hat z},
\beq\label{smooth part limit}
\lim_{k\to\infty} \frac1{2r_k}\int\eta \int_{C_{r_k}^\p} \bigl(d_\e(w+x)-d_\e(w+y)\bigr)\,dz^{r_k}_{(w,t)} = \int\eta\,\nabla d_\e\cdot\widehat z.
\eeq
We now estimate $R_{\e,k}$. We write
\[
\begin{aligned}
R_{\e,k} ={}& \frac1{2r_k}\int\!\int_{C_{r_k}^\p} \Bigl[ (\eta q_{\e,k})(w+x,t) -(\eta q_{\e,k})(w+y,t) \Bigr] \,dz^{r_k}_{(w,t)}(x,y)\,dw\,dt \\
&+ \frac1{2r_k}\int\!\int_{C_{r_k}^\p} \bigl(\eta(w,t)-\eta(w+x,t)\bigr) q_{\e,k}(w+x,t) \,dz^{r_k}_{(w,t)}(x,y)\,dw\,dt \\
&+ \frac1{2r_k}\int\!\int_{C_{r_k}^\p} \bigl(\eta(w+y,t)-\eta(w,t)\bigr) q_{\e,k}(w+y,t) \,dz^{r_k}_{(w,t)}(x,y)\,dw\,dt .
\end{aligned}
\]
By the definition of the nonlocal divergence, the first term is
\[
-\int \eta q_{\e,k}\, d(\mathfrak{Div}_\p^{r_k}z^{r_k}).
\]
Since $\operatorname{supp}\eta$ is compactly contained in the complement of $E$, there exists $\de_0>0$ such that $d\ge2\de_0$ on $\operatorname{supp}\eta$. By the uniform
separation property following Proposition~\ref{distance-function-convergence}, $d_{r_k}\ge\de_0$ on $\operatorname{supp}\eta$ for all sufficiently large $k$. Hence condition
$(\ND)$ gives
\[
\mathfrak{Div}_\p^{r_k}z^{r_k} \le c_{\de_0} \qquad\text{as measures on }\operatorname{supp}\eta,
\]
where $c_{\de_0}$ is the constant appearing in condition $(\ND)$. Since $q_{\e,k}\le0$ on $K\times J$ and $\eta\ge0$, the function $-\eta q_{\e,k}$ is non-negative. Therefore
\[
-\int\eta q_{\e,k}\, d(\mathfrak{Div}_\p^{r_k}z^{r_k}) \le c_{\de_0}\int(-\eta q_{\e,k})\,dw\,dt.
\]
Moreover, by the definition of $m_{\e,k}$,
\[
0\le -q_{\e,k}(\cdot,t) \le \osc_{K} \bigl(d_{r_k}(\cdot,t)-d_\e(\cdot,t)\bigr) \qquad\text{on }K.
\]
The last two terms in the decomposition of $R_{\e,k}$ are estimated using the Lipschitz continuity of $\eta$, the fact that $\po(x),\po(y)\le r_k$ for $(x,y)\in C^\p_{r_k}$, and
$z^{r_k}_{(w,t)}(C_{r_k}^\p)\le1$. Since $K$ contains the $r_k$-neighbourhood of $\operatorname{supp}\eta$ (for $k$ large), we obtain
\beq\label{q error estimate}
R_{\e,k} \le C\int_J \osc_{K} \bigl(d_{r_k}(\cdot,t)-d_\e(\cdot,t)\bigr)\,dt,
\eeq
where $C$ is independent of $k$ and $\e$.

For every $t\notin\mathcal N$, Proposition~\ref{distance-function-convergence} gives $d_{r_k}(\cdot,t)\to d(\cdot,t)$ locally uniformly. Moreover, since $d(\cdot,t)$ is uniformly
Lipschitz with respect to $\pso$,
\[
\|d_\e(\cdot,t)-d(\cdot,t)\|_{L^\infty(K)} \le C\e
\]
uniformly in $t$. Hence
\[
\limsup_{k\to\infty} \osc_{K} \bigl(d_{r_k}(\cdot,t)-d_\e(\cdot,t)\bigr) \le C\e
\]
for a.e. $t\in J$. The quantities on the left are uniformly bounded, by the uniform spatial Lipschitz bounds. Therefore, by dominated convergence in time and
\eqref{q error estimate},
\beq\label{q error limit}
\limsup_{k\to\infty}R_{\e,k}\le C\e.
\eeq
Combining \eqref{smooth part limit} and \eqref{q error limit}, we obtain
\[
\limsup_{k\to\infty}I_k \le \int\eta\,\nabla d_\e\cdot\widehat z+C\e
\]
and, by the arbitrariness of $\e$, we obtain \eqref{CH limsup target}. Combining \eqref{CH liminf start} and \eqref{CH limsup target} gives
\[
\int\eta\,\widehat z\cdot\nabla d \ge \int\eta\,\p(\nabla d)
\]
for every non-negative $\eta$. Thus $\widehat z\cdot\nabla d\ge\p(\nabla d)$ a.e. outside $E$, which concludes the proof.
\end{proof}

\begin{proof}[Proof of Theorem~\ref{asymptotic-behavior-theorem}]
Proposition~\ref{distance-function-convergence} and estimate \eqref{72} give the initial condition and the time-continuity properties required in
Definition~\ref{classic-superflow-def}. Proposition~\ref{divergence-mean-inequality-measures-prop} gives the distributional inequality and the local upper bound on
$(\Div\widehat z)^+$, while Proposition~\ref{Cahn-hoffman-proposition-asymptotic} gives $\widehat z\in\pa\p(\nabla d)$ a.e. away from $E$. Thus all the conditions in
Definition~\ref{classic-superflow-def} are satisfied. Finally, recall that the same argument applies if $(E_r)$ is a family of variational superflows (in the sense of
Definition~\ref{def:variationalflows}) since they satisfy condition $(\ND)$.
\end{proof}

We conclude this section by showing that in the isotropic case condition $(\ND)$ is always satisfied. To this aim we establish the following.

\begin{proposition}[Uniform divergence bound for the Euclidean distance]
\label{isotropic-divergence-bound} Assume that $\p=|\cdot|$ is the Euclidean norm. Let $K\subset\R^N$ be nonempty and closed, $d(x):=\dist(x,K)$, and let $z$ be a generalised
Cahn-Hoffman field for $d$ (with respect to the fixed scale $r$). Then, for every $\de\ge4r$,
\beq\label{isotropic divergence estimate}
\divr_\p z\le \frac{C_N}{\de} \quad\text{on }\{d\ge\de\},
\eeq
where $C_N>0$ depends only on $N$.
\end{proposition}

\begin{proof} Fix $\de>0$ and let $D_\de$ denote the full-measure subset of $\{d\ge\de\}$  of differentiability points of $d$. If $x\in D_\de$ and
$p$ is a projection of $x$ onto $K$, then $d(x+h)\le |x+h-p|$, with equality at $h=0$. Differentiability at $x$ therefore gives \(\nabla d(x)=\frac{x-p}{d(x)}\) and uniqueness of
$p$. Set $\nu(x):=\nabla d(x)$ and consider the map $T(x):=x-r\nu(x)$ on $D_\de$. We now split the proof into two steps.

\noindent\textit{Step 1.} (Semiconcavity estimates)
Consider
\[
f(x):=d^2(x)-|x|^2 =\inf_{p\in K}\bigl(|p|^2-2p\cdot x\bigr).
\]
Since $f$ is the infimum of affine functions, it is concave. Moreover, at every $x\in D_\de$, $\nabla f(x)=2d(x)\nu(x)-2x$. Hence the monotonicity inequality for the gradient of a
concave function gives, for $x,y\in D_\de$,
\[
\begin{aligned}
0 &\ge(\nabla f(x)-\nabla f(y))\cdot(x-y)\\
&=2\bigl(d(x)\nu(x)-d(y)\nu(y)\bigr)\cdot(x-y) -2|x-y|^2,
\end{aligned}
\]
and therefore $\bigl(d(x)\nu(x)-d(y)\nu(y)\bigr)\cdot(x-y)\le |x-y|^2$. Using this inequality, the $1$-Lipschitz continuity of $d$, and $|\nu(y)|=1$, we obtain
\[
\begin{aligned}
d(x)(\nu(x)-\nu(y))\cdot(x-y) &=\bigl(d(x)\nu(x)-d(y)\nu(y)\bigr)\cdot(x-y)\\
&\quad +(d(y)-d(x))\nu(y)\cdot(x-y) \le 2|x-y|^2.
\end{aligned}
\]
Since $d(x)\ge\de$, $(\nu(x)-\nu(y))\cdot(x-y)\le2|x-y|^2/\de$. Thus,
\[
(T(x)-T(y))\cdot(x-y)\ge\big(1-\tfrac{2r}{\de}\big)|x-y|^2, \qquad |T(x)-T(y)|\ge \big(1-\tfrac{2r}{\de}\big)|x-y|.
\]
In particular, $T$ is one-to-one on $D_\de$ and $T^{-1}$ is Lipschitz on $T(D_\de)$. Consequently,
\beq\label{T measure estimates}
|T(A)|\ge\big(1-\tfrac{2r}{\de}\big)^N|A|, \qquad |T^{-1}(B)|\le\big(1-\tfrac{2r}{\de}\big)^{-N}|B|
\eeq
for Borel $A\subset D_\de$ and $B\subset T(D_\de)$.

\smallskip
\noindent\textit{Step 2.} (Estimates of the marginals)
Let $\mu^{\max}$ and $\mu^{\min}$ denote the two marginals of $z$, defined by
\beq\label{marginals}
\begin{aligned}
\int\vphi\,d\mu^{\max} &:=\int_{\R^N}\int_{B_r\times B_r} \vphi(w+x)\,dz_w(x,y)\,dw,\\
\int\vphi\,d\mu^{\min} &:=\int_{\R^N}\int_{B_r\times B_r} \vphi(w+y)\,dz_w(x,y)\,dw,
\end{aligned}
\eeq
for every Borel function $\vphi$. By definition, $\divr_\p z=(\mu^{\min}-\mu^{\max})/(2r)$.

Let $w$ be a differentiability point of $d$ with $d(w)>r$, and let $p$ be its projection onto $K$. Then, $T(w)=p+(d(w)-r)\nu(w)$ and $d(T(w))=d(w)-r$. Moreover, $T(w)$ is the
unique point of minimum of $d$ on $B_r(w)$. Indeed, for $y\in B_r(w)$, $d(y)\ge d(w)-|y-w|\ge d(w)-r$. If equality holds and $q$ is a projection of $y$ onto $K$, then
\[
d(w)\le|w-q|\le|w-y|+|y-q| \le r+d(y)=d(w).
\]
Thus $q$ is also a projection of $w$, hence $q=p$; all the preceding inequalities are equalities, and equality in the triangle inequality forces $y=T(w)$. Since the oscillation of
$d$ on $B_r(w)$ is positive, the defining Cahn-Hoffman identity also gives $z_w(B_r\times B_r)=1$ and implies that the second marginal of $z_w$ is the Dirac mass $\de_{T(w)-w}$;
that is,
\[
\int_{B_r\times B_r} \vphi(w+y)\,dz_w(x,y)=\vphi(T(w))
\]
for any Borel function $\vphi$. In turn, if $A\subset\{d\ge\de\}$ is Borel, recalling that $\de\geq 4r$, we have by \eqref{marginals}
\begin{multline*}
\mu^{\min}(A)=\int_{\R^N}\int_{B_r\times B_r} \chi_A(w+y)\,dz_w(x,y)\,dw=\int_{\{d\geq \frac34\de\}}\int_{B_r\times B_r} \chi_A(w+y)\,dz_w(x,y)\,dw \\
=\int_{\{d\geq \frac34\de\}}\chi_A(T(w))\,dw =|\{w\in D_{\frac34\de}:T(w)\in A\}| \le\big(1-\tfrac{8r}{3\de}\big)^{-N}|A|,
\end{multline*}
where the last inequality follows from \eqref{T measure estimates}.

Conversely, let $w\in T(D_\de)$, set temporarily $x:=T^{-1}(w)$, $s:=d(x)$ and write $p=x-s\nu(x)$ for the projection of $x$ onto $K$. Then $d(w)=s-r>r$, and for every
$y\in B_r(w)$,
\[
d(y)\le|y-p|\le|y-w|+|w-p| \le r+s-r=s=d(x).
\]
If equality holds, then necessarily $p,w,y$ are collinear with $w$ between $p$ and $y$; hence $y=x$. Thus $x=T^{-1}(w)$ is the unique point of maximum of $d$ on $B_r(w)$. Moreover,
arguing as before, we have $z_{w}(B_r\times B_r)=1$ and thus the first marginal of $z_{w}$ is given by the Dirac mass $\de_{T^{-1}(w)-w}$. Let now $A\subset\{d\ge\de\}$ be a Borel
set and define $A_0:=A\cap D_\de$. Then, $|A|=|A_0|$, $T(A_0)\subset T(D_\de)$ and, using again \eqref{marginals}, we have
\begin{multline*}
\mu^{\max}(A)=\int_{\R^N}\int_{B_r\times B_r} \chi_A(w+x)\,dz_w(x,y)\,dw\geq \int_{T(D_\de)}\int_{B_r\times B_r} \chi_{A_0}(w+y)\,dz_w(x,y)\,dw \\
=\int_{T(D_\de)}\chi_{A_0}(T^{-1}(w))\,dw =|T(A_0)\cap T(D_\de)|=|T(A_0)| \geq \big(1-\tfrac{2r}{\de}\big)^N|A|
\end{multline*}
where we used again \eqref{T measure estimates} and the fact that $|A|=|A_0|$. Combining the two marginal estimates gives
\[
\divr_\p z(A) \le\frac{\left(1-\frac{8r}{3\de}\right)^{-N} -\left(1-\frac{2r}{\de}\right)^N}{2r}|A|\leq \frac{\left(1-\frac{8r}{3\de}\right)^{-N}
-\left(1-\frac{8r}{3\de}\right)^N}{2r}|A|
\]
Set
\[
s:=\frac{r}{\de}, \qquad F(s):= \left(1-\frac83s\right)^{-N} -\left(1-\frac83s\right)^N
\]
and note that $s\in[0,1/4]$ (as $\de\ge4r$), $F(0)=0$ and $F$ is Lipschitz on $[0,1/4]$. Denote by $L=L(N)$ its Lipschitz constant. Thus,
\[
\frac{ \left(1-\frac{8r}{3\de}\right)^{-N} -\left(1-\frac{8r}{3\de}\right)^N }{2r} = \frac{F(r/\de)}{2r} \le \frac{L}{2\de}.
\]
Setting $C_N:=L/2$, we have thus shown that $\divr_\p z(A)\le C_N|A|/\de$ for any Borel set contained in $\{d\geq \de\}$. This concludes the proof.
\end{proof}

\begin{remark}\label{isotropic-remark}
Proposition~\ref{isotropic-divergence-bound} shows that in the isotropic case condition $(\ND)$ holds for any family $r\downarrow0$: for every $\de>0$, estimate
\eqref{isotropic divergence estimate} holds uniformly as soon as $r\le\de/4$. Hence Theorem~\ref{asymptotic-behavior-theorem} applies to arbitrary families of isotropic weak
superflows.
\end{remark}


\bibliographystyle{siam}
\bibliography{reference2}

\begin{thebibliography}{10}

\bibitem{ATW}
{\sc F.~Almgren, J.~E. Taylor, and L.~Wang}, {\em Curvature-driven flows: a
  variational approach}, SIAM Journal on Control and Optimization, 31 (1993),
  pp.~387--438.

\bibitem{AFP}
{\sc L.~Ambrosio, N.~Fusco, and D.~Pallara}, {\em Functions of Bounded
  Variation and Free Discontinuity Problems}, Oxford Mathematical Monographs,
  Oxford University Press, Oxford, 2000.

\bibitem{BKLMP}
{\sc M.~Barchiesi, S.~H. Kang, T.~M. Le, M.~Morini, and M.~Ponsiglione}, {\em A
  variational model for infinite perimeter segmentations based on lipschitz
  level set functions: denoising while keeping finely oscillatory boundaries},
  Multiscale Modeling \& Simulation, 8 (2010), pp.~1715--1741.

\bibitem{BLS24}
{\sc L.~Bungert, T.~Laux, and K.~Stinson}, {\em A mean curvature flow arising
  in adversarial training}, Journal de Math\'ematiques Pures et Appliqu\'ees,
  192 (2024), p.~103625.

\bibitem{BS26}
{\sc L.~Bungert and L.~Schmitt}, {\em Duality for the adversarial total
  variation}, 2026.

\bibitem{BS24}
{\sc L.~Bungert and K.~Stinson}, {\em Gamma-convergence of a nonlocal perimeter
  arising in adversarial machine learning}, Calculus of Variations and Partial
  Differential Equations, 63 (2024), p.~114.

\bibitem{CMM08}
{\sc F.~Cagnetti, M.~G. Mora, and M.~Morini}, {\em A second order minimality
  condition for the mumford-shah functional}, Calculus of Variations and
  Partial Differential Equations, 33 (2008), pp.~37--74.

\bibitem{CMR26}
{\sc F.~Cagnetti, M.~Morini, and D.~Reggiani}, {\em A distributional approach
  to nonlocal curvature flows}, Communications in Partial Differential
  Equations, 51 (2026), pp.~364--413.

\bibitem{CDNP21}
{\sc A.~Cesaroni, L.~De~Luca, M.~Novaga, and M.~Ponsiglione}, {\em Stability
  results for nonlocal geometric evolutions and limit cases for fractional mean
  curvature flows}, Communications in Partial Differential Equations, 46
  (2021), pp.~1344--1371.

\bibitem{CN17}
{\sc A.~Cesaroni and M.~Novaga}, {\em Isoperimetric problems for a nonlocal
  perimeter of minkowski type}, Geometric Flows, 2 (2017), pp.~86--93.

\bibitem{Chambolle}
{\sc A.~Chambolle}, {\em An algorithm for mean curvature motion}, Interfaces
  and Free Boundaries, 6 (2004), pp.~195--218.

\bibitem{CDGM26}
{\sc A.~Chambolle, D.~De~Gennaro, and M.~Morini}, {\em Discrete-to-continuum
  crystalline curvature flows}, Analysis \& PDE, 19 (2026), pp.~783--822.

\bibitem{CLL14}
{\sc A.~Chambolle, S.~Lisini, and L.~Lussardi}, {\em A remark on the
  anisotropic outer minkowski content}, Advances in Calculus of Variations, 7
  (2014), pp.~241--266.

\bibitem{CMNP19}
{\sc A.~Chambolle, M.~Morini, M.~Novaga, and M.~Ponsiglione}, {\em Existence
  and uniqueness for anisotropic and crystalline mean curvature flows}, Journal
  of the American Mathematical Society, 32 (2019), pp.~779--824.

\bibitem{CMP12}
{\sc A.~Chambolle, M.~Morini, and M.~Ponsiglione}, {\em A nonlocal mean
  curvature flow and its semi-implicit time-discrete approximation}, SIAM
  Journal on Mathematical Analysis, 44 (2012), pp.~4048--4077.

\bibitem{CMP15}
\leavevmode\vrule height 2pt depth -1.6pt width 23pt, {\em Nonlocal curvature
  flows}, Archive for Rational Mechanics and Analysis, 218 (2015),
  pp.~1263--1329.

\bibitem{CMP17}
\leavevmode\vrule height 2pt depth -1.6pt width 23pt, {\em Existence and
  uniqueness for a crystalline mean curvature flow}, Communications on Pure and
  Applied Mathematics, 70 (2017), pp.~1084--1114.

\bibitem{FonsecaLeoniBook}
{\sc I.~Fonseca and G.~Leoni}, {\em Modern Methods in the Calculus of
  Variations: $L^p$ Spaces}, Springer Monographs in Mathematics, Springer, New
  York, 2007.

\bibitem{GP16}
{\sc Y.~Giga and N.~Po{\v{z}}{\'a}r}, {\em A level set crystalline mean
  curvature flow of surfaces}, Advances in Differential Equations, 21 (2016),
  pp.~631--698.

\bibitem{GP18}
\leavevmode\vrule height 2pt depth -1.6pt width 23pt, {\em Approximation of
  general facets by regular facets with respect to anisotropic total variation
  energies and its application to crystalline mean curvature flow},
  Communications on Pure and Applied Mathematics, 71 (2018), pp.~1461--1491.

\bibitem{GP20}
\leavevmode\vrule height 2pt depth -1.6pt width 23pt, {\em Viscosity solutions
  for the crystalline mean curvature flow with a nonuniform driving force
  term}, SN Partial Differential Equations and Applications, 1 (2020), p.~39.

\bibitem{GP22}
\leavevmode\vrule height 2pt depth -1.6pt width 23pt, {\em Motion by
  crystalline-like mean curvature: a survey}, Bulletin of Mathematical
  Sciences, 12 (2022), p.~2230004.

\bibitem{Rockafellar}
{\sc R.~T. Rockafellar}, {\em Convex Analysis}, vol.~28 of Princeton
  Mathematical Series, Princeton University Press, Princeton, NJ, 1970.

\bibitem{S93}
{\sc H.~M. Soner}, {\em Motion of a set by the curvature of its boundary},
  Journal of Differential Equations, 101 (1993), pp.~313--372.

\end{thebibliography}

\end{document}